\documentclass{amsart}
\usepackage{ifthen, citesort}
\usepackage{amssymb}

\DeclareMathOperator{\tr}{tr}
\DeclareMathOperator{\divergence}{div}
\DeclareMathOperator{\divergenz}{div}
\DeclareMathOperator{\sgn}{sgn}
\DeclareMathOperator{\diag}{diag}
\def\R{\mathbb{R}}
\def\N{\mathbb{N}}
\def\S{\mathbb{S}}
\def\dt{\frac{d}{dt}}
\def\fracp#1#2{\frac{\partial #1}{\partial #2}}

\def\theta{\vartheta}
\def\phi{\varphi}
\def\epsilon{\varepsilon}
\newcommand{\A}[1]{\ifthenelse{#1 = 2}{\lvert A\rvert^{#1}}{\tr A^{#1}}}
\newcommand{\Akl}[2]{\ifthenelse{#1 = 2}%
{\left(\lvert A\rvert^{#1}\right)^{#2}}%
{\left(\tr A^{#1}\right)^{#2}}}
\def\B#1{\tr A^{#1}}
\def\hij{h_{ij;\,k}^2}
\def\hii{h_{ii;\,k}h_{jj;\,k}}
\def\li{\lambda_i}
\def\lj{\lambda_j}
\def\lk{\lambda_k}
\def\lx{\lambda_1}
\def\ly{\lambda_2}
\def\opr#1{\dt{#1}-F^{ij}{\left(#1\right)}_{;\,ij}}
\def\oprokl#1{\dt{#1}-F^{ij}{#1}_{;\,ij}}
\def\fracd#1#2{\displaystyle\frac{\displaystyle #1}{\displaystyle #2}}
\def\fracm#1#2{\textstyle\frac{\scriptstyle #1}{\scriptstyle #2}}
\long\def\etamatrix#1#2#3{\begin{pmatrix}\eta_{11}\\
\eta_{22}\end{pmatrix}^{\text{tr}}\begin{pmatrix}#1 &#2\\
#2& #3\end{pmatrix}\begin{pmatrix}\eta_{11}\\
\eta_{22}\end{pmatrix}}
\long\def\symmatrix#1#2#3{\begin{pmatrix}#1 &#2\\
#2& #3\end{pmatrix}}
\long\def\umbruch{{\displaybreak[1]}}

\newtheorem{theorem}{Theorem}[section]
\newtheorem{lemma}[theorem]{Lemma}
\newtheorem{proposition}[theorem]{Proposition}
\newtheorem{corollary}[theorem]{Corollary}

\theoremstyle{definition}

\theoremstyle{remark}

\numberwithin{equation}{section}
\newcommand{\abs}[1]{\left\lvert#1\right\rvert}

\begin{document}

\title{Surfaces contracting with speed $|A|^2$}

\author{Oliver C. Schn\"urer}
\address{FU Berlin, Arnimallee 2-6, 14195 Berlin, Germany}
\curraddr{}
\email{Oliver.Schnuerer@math.fu-berlin.de}
\thanks{For a short version of this paper see
  J. Differential Geom. \textbf{71} (2005), no.~3, 347--363.}


\date{September 2004, revised February 2006.}

\dedicatory{}


\begin{abstract}
We show that strictly convex surfaces contracting with normal 
velocity equal to $\A2$ shrink to a point in finite time.
After appropriate rescaling, they converge to spheres.
We describe our algorithm to find the main test function.
\end{abstract}

\maketitle

\tableofcontents

\section{Introduction}

We consider closed strictly convex surfaces $M_t$ in $\R^3$ 
that contract with normal velocity equal to the square
of the norm of the second fundamental form
\begin{equation}\label{flow eqn}
\dt X=-\A2\nu.
\end{equation}
This is a parabolic flow equation. We obtain a solution
on a maximal time interval $[0,\,T)$, $0<T<\infty$. 
For $t\uparrow T$, the surfaces converge to a point. After
appropriate rescaling, they converge to a round sphere.
We say that the surfaces $M_t$ converge to a ``round point''.
The key step in the proof, Theorem \ref{A2 mon thm}, is to show that
\begin{equation}\label{mon groe}
\max\limits_{M_t}\left(\frac{(\lx+\ly)(\lx-\ly)^2}{\lx\ly}\right)
\end{equation}
is non-increasing in time. \par 
Here, we used standard notation 
as explained in Section \ref{nota sec}.

Our main theorem is
\begin{theorem}\label{main thm}
For any smooth closed strictly convex surface $M$ in $\R^3$, there
exists a smooth family of surfaces $M_t$, $t\in[0,T)$,
solving \eqref{flow eqn} with $M_0=M$. For $t\uparrow T$,
$M_t$ converges to a point $Q$. The rescaled surfaces
$(M_t-Q)\cdot(6(T-t))^{-1/3}$ converge smoothly to the
unit sphere $\S^2$.
\end{theorem}

We will also consider other normal velocities for which
similar results hold. Therefore, we have to find 
quantities like \eqref{mon groe} that are monotone
during the flow and vanish precisely for spheres.
In general, this is a complicated issue. 
In order to find these test quantities, we used an algorithm
that checks, based on randomized tests, whether possible
candidates fulfill certain inequalities. These inequalities
guarantee especially that we can apply the maximum principle
to prove monotonicity. We used that algorithm only to propose
useful quantities. The presented proofs do not depend on it. 
So far, all candidates turned out to be
appropriate for proving convergence to a round point. 
Our algorithm yields also candidates for many other
normal velocities. We have only included a discussion of
some interesting normal velocities. Moreover, for a
fixed normal velocity, there are mostly several candidates
for monotone quantities. In these cases, we have picked those
involving not too complicated polynomials of low homogeneity. 

In Table \ref{table}, we have collected some normal
velocities $F$ ($1^{\text{st}}$ column) and quantities
$w$ ($2^{\text{nd}}$ column) such that $\max_{M_t}w$ is
non-increasing in time for surfaces contracting
with normal velocity $F$. In each case, we obtain
convergence to round points for smooth closed strictly 
convex initial surfaces $M_0$. 

\begin{table}
\def\platz{\raisebox{0em}[2.2em][1.5em]{\rule{0em}{2em}}}
$$
\begin{array}{|c||c|}\hline
\A2 & \platz\fracd{(\lx+\ly)(\lx-\ly)^2}{\lx\ly}\\\hline
K {\text{ \cite{AndrewsStones}}} & \platz(\lx-\ly)^2\\\hline
H^2 & \platz\fracd{(\lx+\ly)^3(\lx-\ly)^2}
{\left(\lx^2+\ly^2\right)\lx\ly}\\\hline
H^3 & \platz\fracd{\left(\lx^2+\lx\ly+\ly^2\right)(\lx+\ly)^2(\lx-\ly)^2}
{\left(\lx^2-\lx\ly+\ly^2\right)\lx\ly}\\\hline
H^4 & \platz\fracd{\left(\lx^2+\lx\ly+\ly^2\right)(\lx+\ly)^6(\lx-\ly)^2}
{\lx^2\ly^2}\\\hline
\begin{minipage}{1.7cm}
\centerline{$\A2+\beta H^2,$}
\centerline{$0\le\beta\le5$} 
\end{minipage}
& \platz\fracd{(\lx+\ly)(\lx-\ly)^2}{\lx\ly}\\\hline
\A3 & \platz\fracd{\left(3\lx^2+2\lx\ly+3\ly^2\right)(\lx-\ly)^2}
{\lx\ly}\\\hline
\begin{minipage}{2.1cm}
\centerline{$\B{\alpha},$}
\centerline{$\alpha=2,\,4,\,5,\,6$}
\end{minipage}
& \platz\fracd{\left(\lx^{\alpha-2}+\ly^{\alpha-2}\right)
(\lx+\ly)(\lx-\ly)^2}
{\lx\ly}\\\hline
H\A2 & \platz\fracd{(\lx+\ly)^2(\lx-\ly)^2}{\lx\ly}\\\hline
\abs A^4 & \platz\fracd{\left(\lx^4+2\lx^3\ly+4\lx^2\ly^2
+2\lx\ly^3+\ly^4\right)(\lx-\ly)^2}{(\lx+\ly)\lx\ly}\\\hline
\end{array}$$
\caption{Monotone quantities}\label{table}
\end{table}

In \cite{HuiskenRoundSphere,Andrews2dnonconcave,%
AndrewsPinching,AndrewsContractCalc}, Gerhard Huisken
and Ben Andrews proved that convex hypersurfaces
contracting with certain normal velocities homogeneous
of degree one converge to ``round points'', i.\,e., they
converge to a point and, after appropriate rescaling,
to round spheres. For homogeneities larger than one,
this was shown by Ben Andrews and Felix Schulze
\cite{Andrews2dnonconcave,OSFelixH2}, if the
initial hypersurfaces are pinched appropriately.
Kaising Tso proved that Gau{\ss} curvature flow
shrinks strictly convex hypersurfaces to points \cite{TsoPoint}.
If the homogeneity is less than one, there are
examples by Koichi Anada, Masayoshi Tsutsumi, and
Ben Andrews, where hypersurfaces do not become
spherical \cite{AndrewsPacific,AnadaCalcHarm,%
AnadaStability}. Expanding flows were studied 
by Claus Gerhardt, John Urbas, Bennett Chow, Dong-Ho Tsai,
Nina Ivochkina, Thomas Nehring, Friedrich Tomi, Knut Smoczyk,
Gerhard Huisken, and Tom Ilmanen \cite{CGFlowSpheres,%
HuiskenIlmanenPenrose,SmoczykAsian2000,UrbasExpandJDG,%
UrbasExpandMZ,KnutHarmonicMCF,ChowAsian,IvochkinaNehringTomi}. 
Similar problems were also studied in manifolds 
(e.\,g.\ \cite{AndrewsBrisbane,%
HuiskenSphereMf,AndrewsContrRiem}) and for anisotropic flow
equations (e.\,g.\ \cite{AndrewsMonotone}). 
It is often required that the normal velocity is a concave function
of the second fundamental form.
There are many papers, concerned with contracting curves, e.\,g.\ by
Michael Gage, Richard Hamilton, Matthew Grayson, and
Steven Altschuler \cite{GraysonRoundPoints,GageHamilton,%
SteveThroughSing}. 

In \cite{AndrewsStones}, Ben Andrews shows that convex surfaces
moving by Gau{\ss} curvature converge to round points. This normal
velocity is homogeneous of degree two in the principal curvatures.
He does not require any pinching condition for the initial 
surface. Our paper extends this result to other flow equations.
We consider also normal velocities of degree larger than one and
do not have to impose any pinching condition on the initial surface.
Any smooth strictly convex surface converges to a round point.

The rest of this paper is organized as follows. In Section 
\ref{nota sec}, we explain our notation. Section \ref{A2 sec}
concerns the proof for the normal velocity $\A2$. We 
describe our algorithm to find test quantities 
in Section \ref{algor sec}. In the
remaining sections, we prove convergence for some other normal
velocities and discuss the expected convergence rate. 

The author wants to thank Shing-Tung Yau at Harvard,
the Alexander von Humboldt foundation, J\"urgen Jost at the
Max Planck Institute in Leipzig and Klaus Ecker at the Free
University Berlin for discussions and support. We also want to
thank John Stalker for telling us about Sturm's theorem,
Felix Schulze for discussing the convergence proof of
\cite{AndrewsStones} and 
Olaf Schn\"urer for pointing out an appropriate basis.
Kashif Rasul told us useful C-compiler options. 

\section{Notation}\label{nota sec}
We use $X=X(x,\,t)$ to denote the embedding vector of a manifold 
$M_t$ into $\R^3$ and $\dt X=\dot X$ for its total time derivative. 
It is convenient to identify $M_t$ and its embedding in $\R^3$.
We choose $\nu$ to be the outer unit normal vector to $M_t$. 
The embedding induces a metric $(g_{ij})$ and
a second fundamental form $(h_{ij})$. We use the Einstein summation
convention. Indices are raised and lowered with respect to the metric
or its inverse $\left(g^{ij}\right)$. The inverse of the second fundamental
form is denoted by $\left(\tilde h^{ij}\right)$. The principal
curvatures $\lx,\,\ly$ are the eigenvalues of the second fundamental
form with respect to the induced metric. A surface is called strictly convex,
if all principal curvatures are strictly positive.
We will assume this throughout the paper.

Symmetric functions of the principal
curvatures are well-defined, we will use the mean curvature
$H=\lx+\ly$, the square of the norm of the second fundamental form
$\A2=\lx^2+\ly^2$, $\B k=\lx^k+\ly^k$, and the Gau{\ss} curvature
$K=\lx\ly$. We write indices, preceded by semi-colons, e.\,g.\ $h_{ij;\,k}$, 
to indicate covariant differentiation with respect to the induced
metric. It is often convenient to choose coordinate systems such
that the metric tensor equals the Kronecker delta, $g_{ij}=\delta_{ij}$,
and $(h_{ij})$ is diagonal, $(h_{ij})=\diag(\lx,\,\ly)$, e.\,g.
$$\sum\lk\hij=\sum\limits_{i,\,j,\,k=1}^2\lk\hij
=h^{kl}h^i_{j;\,k}h^j_{i;\,l}
=h_{rs}h_{ij;\,k}h_{ab;\,l}g^{ia}g^{jb}g^{rk}g^{sl}.$$
Whenever we use this notation, we will also assume that we have 
fixed such a coordinate system. We will only use Euclidean coordinate
systems for $\R^3$ so that $h_{ij;\,k}$ is symmetric according to
the Codazzi equations.

A normal velocity $F$ can be considered as a function of $(\lx,\,\ly)$
or $(h_{ij},\,g_{ij})$. We set $F^{ij}=\fracp{F}{h_{ij}}$,
$F^{ij,\,kl}=\fracp{^2F}{h_{ij}\partial h_{kl}}$. 
Note that in coordinate
systems with diagonal $h_{ij}$ and $g_{ij}=\delta_{ij}$ as mentioned
above, $F^{ij}$ is diagonal. 
For $F=\A2$, we have $F^{ij}=2h^{ij}=2\li g^{ij}$. 

Recall, see e.\,g.\ \cite{HuiskenRoundSphere,OSKnutCrelle,OSMM}, 
that for a hypersurface moving according to $\dt X=-F\nu$, we have
\begin{align}
\label{g evol}
\dt g_{ij}=&-2Fh_{ij},\displaybreak[1]\\
\label{h evol}
\dt h_{ij}=&F_{;\,ij}-Fh_i^kh_{kj},\displaybreak[1]\\
\dt\nu^\alpha=&g^{ij}F_{;\,i}X^\alpha_{;\,j},
\end{align}
where Greek indices refer to components in the ambient space $\R^3$.
In order to compute evolution equations, we use the
Gau{\ss} equation and the Ricci identity for the second fundamental
form
\begin{align}
\label{Riem}
R_{ijkl}=&h_{ik}h_{jl}-h_{il}h_{jk},\displaybreak[1]\\
\label{Ricci}
h_{ik;\,lj}=&h_{ik;\,jl}+h^a_kR_{ailj}+h^a_iR_{aklj}.
\end{align}
We will also employ the Gau{\ss} formula and the Weingarten
equation
$$X^\alpha_{;\,ij}=-h_{ij}\nu^\alpha\qquad\text{and}\qquad
\nu^\alpha_{;\,i}=h^k_iX^\alpha_{;\,k}.$$

For tensors $A$ and $B$, $A_{ij}\ge B_{ij}$ means that 
$(A_{ij}-B_{ij})$ is positive definite.
Finally, we use $c$ to denote universal, estimated constants. 

\section{Surfaces Flowing With Speed $|A|^2$}\label{A2 sec}

\subsection{Convergence to a Point}
It is known, that \eqref{flow eqn} is a parabolic evolution equation
for strictly convex initial data and that is has a solution on a maximal
time interval. 

We show that $M_t$ stays uniformly strictly convex.
The following lemma is similar to results in \cite{AndrewsContractCalc}.
\begin{lemma}\label{A2 conv pre}
For a smooth closed strictly convex surface $M$ in $\R^3$, 
flowing according to $\dot X=-\A2\nu$, the minimum of the
principal curvatures is non-decreasing.
\end{lemma}
\begin{proof}
Consider $M_{ij}=h_{ij}-\epsilon g_{ij}$ with $\epsilon>0$ so small
that $M_{ij}$ is positive semi-definite for some time $t_0$. We wish
to show that $M_{ij}$ is positive semi-definite for $t>t_0$. Combine
\eqref{h evol}, \eqref{Riem}, and \eqref{Ricci} to obtain
$$\dt h_{ij}-F^{kl}h_{ij;\,kl}=2\A3 h_{ij}-3\A2h^k_ih_{kj}
+2g^{kr}g^{ls}h_{kl;\,i}h_{rs;\,j}.$$
In the evolution equation for $M_{ij}$, we drop the positive definite 
terms involving derivatives of the second fundamental form
$$\dt M_{ij}-F^{kl}M_{ij;\,kl}\ge 2\A3 h_{ij}-3\A2h^k_ih_{kj}
+2\epsilon\A2 h_{ij}.$$
Let $\xi$ be a zero eigenvalue of $M_{ij}$ with $\abs{\xi}=1$,
$M_{ij}\xi^j=h_{ij}\xi^j-\epsilon g_{ij}\xi^j=0$.
So we obtain in a point with $M_{ij}\ge0$
\begin{align*}
\left(2\A3h_{ij}-3\A2 h^k_ih_{kj}+2\epsilon\A2 h_{ij}\right)
\xi^i\xi^j=&2\epsilon\A3-3\epsilon^2\A2+2\epsilon^2\A2\\
=&2\epsilon\A3-\epsilon^2\A2\\
\ge&2\epsilon^2\A2-\epsilon^2\A2>0
\end{align*}
and the maximum principle for tensors 
\cite{ChowLuMP,HamiltonThree,HamiltonFour86}
gives the result.
\end{proof}

The next result shows that $\A2$ stays uniformly bounded as
long as $M_t$ encloses a ball of fixed positive radius. A similar
estimate is used in \cite{TsoPoint}.
\begin{lemma}\label{vel bound}
For a strictly convex solution of \eqref{flow eqn}, $\A2$ is
uniformly bounded in terms of the radius $R$ of an enclosed 
sphere $B_{R}(x_0)$,
$\max_{M_0}\,\frac{\A2}{\langle X-x_0,\,\nu\rangle-\frac12R}$,
and $\max_{M_0}|X-x_0|$. More precisely, we have
\begin{equation}\label{A2 bound}
\sup\limits_{t}\max\limits_{M_t}\,\A2\le\max\left\{
\max_{M_0}|X-x_0|\cdot\max_{M_0}\frac{\A2}{\langle X-x_0,\,\nu\rangle
-\frac12R},\,\frac{18}{R^2}\right\}.
\end{equation}
\end{lemma}
\begin{proof}
We may assume that $x_0=0$. Let $\alpha=\frac12R$. Then
$\alpha$ is a positive lower bound for $\langle X,\,\nu\rangle-\alpha$.
Standard computations \cite{HuiskenRoundSphere,OSMM,OSKnutCrelle} 
yield the evolution equations
\begin{align*}
\dt X^\beta-F^{ij}X^\beta_{;\,ij}=&\A2\nu^\beta,\displaybreak[1]\\
\dt\nu^\beta-F^{ij}\nu^\beta_{;\,ij}=&2\A3\nu^\beta,\displaybreak[1]\\
\dt\langle X,\,\nu\rangle-F^{ij}\langle X,\,\nu\rangle_{;\,ij}=&
-3\A2+2\A3\langle X,\,\nu\rangle,\displaybreak[1]\\
\opr{\A2}=&2\A2\A3.
\end{align*}
In a critical point of $\frac{\A2}{\lvert X,\,\nu\rvert-\alpha}$,
we obtain
$$\dt\log\frac{\A2}{\langle X,\,\nu\rangle-\alpha}
-F^{ij}\left(\log\frac{\A2}{\langle X,\,\nu\rangle-\alpha}\right)_{;\,ij}=
\frac1{\langle X,\,\nu\rangle-\alpha}\left(3\A2-2\A3\alpha\right).$$
Note that $\langle X,\,\nu\rangle-\alpha\le\max_{M_0}
\abs X$ as a sphere of radius $\max_{M_0}\abs X$, 
centered at the origin, will enclose
any $M_t$.
We only have to prove that we preserve the bound in 
Equation \eqref{A2 bound}, when 
$\max_{M_t}\frac{\A2}{\langle X,\,\nu\rangle-\alpha}$ increases.
Then we have $0\le3\A2-2\A3\alpha$ at a point, where
$\max_{M_t}\frac{\A2}{\langle X,\,\nu\rangle-\alpha}$ 
is attained. This inequality and elementary calculations 
for convex surfaces give 
$$\A2\le2^{1/3}\cdot\Akl3{2/3}\le
2\frac{\Akl32}{\Akl22}\le\frac9{2\alpha^2}$$
at such a maximum point and the Lemma follows. 
\end{proof}

We obtain that the second fundamental form of the surface stays bounded
as long as $M_t$ encloses some ball. The estimates of Krylov, Safonov,
Evans (see also \cite{Andrews2dKrylov}), 
and Schauder imply that the solution stays smooth. Then,
similarly as in \cite{TsoPoint}, 
the positive lower bound on the minimum principal curvature
implies that the surfaces converge to a point in finite time.

\subsection{A Monotone Quantity}\label{mon sec}
\begin{theorem}\label{A2 mon thm}
For a family of smooth closed strictly convex surfaces $M_t$ 
in $\R^3$ flowing according to $\dot X=-\A2\nu$, 
\begin{equation}\label{A2 quantity}
\max\limits_{M_t}\frac{(\lambda_1+\lambda_2)(\lambda_1-\lambda_2)^2}
{2\lambda_1\lambda2}
=\max\limits_{M_t}\frac{H\cdot\left(2\A2-H^2\right)}{H^2-\A2}
\equiv\max\limits_{M_t}w
\end{equation}
is non-increasing in time.
\end{theorem}

An immediate consequence of this theorem is
\begin{corollary}
The only homothetically shrinking smooth closed strictly convex
surfaces $M_t$, solving the flow equation $\dot X=-\A2\nu$ in 
$\R^3$, are spheres.
\end{corollary}
\begin{proof}
The quantity $\frac{(\lx+\ly)(\lx-\ly)^2}{\lx\ly}$ 
is positive homogeneous of degree one
in the principal curvatures and non-negative. If $M$ is 
homothetically shrinking, Theorem \ref{A2 mon thm} implies that
$(\lx+\ly)(\lx-\ly)^2=0$ everywhere. Thus $M_t$ is umbilic and 
\cite[Lemma 7.1]{Spivak4} implies that $M_t$ is a sphere.
\end{proof}

\begin{proof}[Proof of Theorem \ref{A2 mon thm}]
We combine \eqref{g evol}, \eqref{h evol}, \eqref{Riem},
and \eqref{Ricci} in order to get the 
following general evolution equation
\begin{align}\label{gen ev eq}
\begin{split}
\dt\B\alpha-F^{kl}\left({\B\alpha}\right)_{;\,kl}=&
\alpha\sum\limits_i F^{ii}\lambda_i^2\B\alpha
+\alpha\left(F-\sum\limits_i F^{ii}\lambda_i\right)
\B{\alpha+1}\\
&-\alpha\sum\limits_{r=0}^{\alpha-2}\sum\limits_{i,\,j,\,k} 
F^{kk}\lambda_i^{\alpha-2-r}
\lambda_j^r\hij\\
&+\alpha\sum\limits_{k,\,l,\,r,\,s,\,i} 
F^{kl,\,rs}h_{kl;\,i}h_{rs;\,i}\lambda_i^{\alpha-1}
\end{split}
\end{align}
for solutions to the flow equation
$$\dt X=-F\nu.$$
Using \eqref{gen ev eq} for $F=\A2$ yields
\begin{align}
\begin{split}\label{H evol}
\oprokl H=&-\left(\A2\right)^2+2H\A3+2\sum\hij
\end{split}
\intertext{and}
\begin{split}\label{A2 evol}
\opr{\A2}=&2\A2\A3.
\end{split}
\end{align}
For the reader's convenience, we include the details to obtain
\eqref{H evol}.
\begin{align*}
H=&g^{ij}h_{ij},&&\umbruch\\
F=&\A2=g^{ij}h_{jk}g^{kl}h_{li},&&\umbruch\\
F^{ij}=&2h^{ij},&&\umbruch\\
\dt H=&-g^{ia}g^{bj}h_{ij}\dt g_{ab}+g^{ij}\dt h_{ij}&&\umbruch\\
=&-g^{ia}g^{bj}h_{ij}\left(-2\A2 h_{ab}\right)&&\\
&+g^{ij}\left(\left(\A2\right)_{;\,ij}-\A2 h^k_ih_{kj}\right)&&
\text{by \eqref{g evol} and \eqref{h evol}}\umbruch\\
=&\Akl22+g^{ij}\left(2h^{kl}h_{kl;\,ij}
+2g^{kr}g^{ls}h_{kl;\,i}h_{rs;\,j}\right)&&\umbruch\\
=&\Akl22+2g^{ij}h^{kl}h_{kl;\,ij}+2\sum\hij&&\umbruch\\
=&\Akl22+2g^{ij}h^{kl}\left(h_{ij;\,kl}
+h^a_kR_{ailj}+h^a_iR_{aklj}\right)&&\\
&+2\sum\hij&&\text{by \eqref{Ricci} and Codazzi}\umbruch\\
=&F^{kl}\left(g^{ij}h_{ij}\right)_{;\,kl}+\Akl22&&\\
&+2g^{ij}h^{kl}h^a_kh_{al}h_{ij}-2g^{ij}h^{kl}h^a_kh_{aj}h_{il}&&\\
&+2g^{ij}h^{kl}h^a_ih_{al}h_{kj}-2g^{ij}h^{kl}h^a_ih_{aj}h_{kl}
+2\sum\hij&&\text{by \eqref{Riem}}\umbruch\\
=&F^{ij}H_{;\,ij}-\Akl22+2H\A3+2\sum\hij.&&
\end{align*}
For the rest of the proof, we consider a critical point of 
$\left.w\right|_{M_t}$ for some $t>0$, where $w>0$.
It suffices to show that $\tilde w:=\log w$ is non-increasing
in such a point. Then our theorem follows. 

We rewrite $\tilde w$ 
\begin{align*}
\tilde w=&\log H+\log\left(2\A2-H^2\right)-\log(H^2-\A2)\\
\equiv&\log A+\log B-\log C.
\end{align*}
In a critical point of $\tilde w$, we obtain
\begin{align*}
\begin{split}
\oprokl{\tilde w}=&\frac1A\left(\oprokl A\right)
+\frac1B\left(\oprokl B\right)\\
&-\frac1C\left(\oprokl C\right)
-\frac1{AB}F^{ij}(A_{;\,i}B_{;\,j}+A_{;\,j}B_{;\,i})\\
\end{split}\\
\intertext{and}
\begin{split}
0=&\frac1HH_{;\,k}+\frac1{2\A2-H^2}\left(2\A2-H^2\right)_{;\,k}
-\frac1{H^2-\A2}\left(H^2-\A2\right)_{;\,k}\\
=&\frac{2\lx^2+\lx\ly+\ly^2}{\lx\left(\lx^2-\ly^2\right)}h_{11;\,k}
+\frac{2\ly^2+\lx\ly+\lx^2}{\ly\left(\ly^2-\lx^2\right)}h_{22;\,k}.
\end{split}
\end{align*}
So we deduce that
$$h_{22;\,1}=\frac{\ly}{\lx}\frac{2\lx^2+\lx\ly+\ly^2}{2\ly^2+\lx\ly+\lx^2}
h_{11;\,1}\equiv a_1h_{11;\,1}$$
and a similar formula holds for $h_{11;\,2}$, 
$h_{11;\,2}=a_2\cdot h_{22;\,2}$. 
We now combine all these results and obtain in a straightforward
calculation
\begin{align*}
\oprokl{\tilde w}=&
\left(\frac1H-\frac{2H}{2\A2-H^2}-\frac{2H}{H^2-\A2}\right)
\cdot\left(\oprokl H\right)\displaybreak[1]\\
&+\left(\frac{2}{2\A2-H^2}+\frac1{H^2-\A2}\right)
\cdot\left(\opr{\A2}\right)\displaybreak[1]\\
&+\left(\frac{6}{2\A2-H^2}+\frac{2}{H^2-\A2}\right)F^{ij}H_{;\,i}H_{;\,j}
\displaybreak[1]\\
&-\frac2{H\cdot\left(2\A2-H^2\right)}
F^{ij}\left(\left(\A2\right)_{;\,i}H_{;\,j}
+\left(\A2\right)_{;\,j}H_{;\,i}\right)
\displaybreak[1]\\
=&-\frac{\lambda_1^4+\lambda_1^3\lambda_2+4\lambda_1^2\lambda_2^2
+\lambda_1\lambda_2^3+\lambda_2^4}{(\lambda_1+\lambda_2)(\lambda_1-
\lambda_2)^2\lambda_1\lambda_2}\cdot\\
&\quad\cdot\left(-\left(\A2\right)^2+2H\A3\right)\displaybreak[1]\\
&+\frac{(\lambda_1+\lambda_2)^2}{2(\lambda_1-\lambda_2)^2
\lambda_1\lambda_2}\cdot
2\A2\A3\displaybreak[1]\\
&-2\frac{\lambda_1^4+\lambda_1^3\lambda_2+4\lambda_1^2\lambda_2^2
+\lambda_1\lambda_2^3+\lambda_2^4}{(\lambda_1+\lambda_2)
(\lambda_1-\lambda_2)^2\lambda_1\lambda_2}
\sum\hij\displaybreak[1]\\
&+2\frac{\lambda_1^2+4\lambda_1\lambda_2+\lambda_2^2}
{(\lambda_1-\lambda_2)^2\lambda_1\lambda_2}
\sum\lk\hii\displaybreak[1]\\
&-\frac{8}{(\lambda_1-\lambda_2)^2(\lambda_1+\lambda_2)}
\sum\lk(\li+\lj)\hii\displaybreak[1]\\
=&-\frac{\lambda_1^8+3\lambda_1^7\lambda_2+4\lambda_1^6\lambda_2^2
+9\lambda_1^5\lambda_2^3-2\lambda_1^4\lambda_2^4}
{(\lambda_1+\lambda_2)(\lambda_1-\lambda_2)^2\lambda_1\lambda_2}
\displaybreak[1]\\
&-\frac{9\lambda_1^3\lambda_2^5
+4\lambda_1^2\lambda_2^6+3\lambda_1\lambda_2^7+\lambda_2^8}
{(\lambda_1+\lambda_2)(\lambda_1-\lambda_2)^2\lambda_1\lambda_2}
\displaybreak[1]\\
&+\frac{\lambda_1^7+2\lambda_1^6\lambda_2+2\lambda_1^5\lambda_2^2
+3\lambda_1^4\lambda_2^3+3\lambda_1^3\lambda_2^4+2\lambda_1^2\lambda_2^5
+2\lambda_1\lambda_2^6+\lambda_2^7}
{(\lambda_1-\lambda_2)^2\lambda_1\lambda_2}\displaybreak[1]\\
&-2\frac{\lx^4+\lx^3\ly+4\lx^2\ly^2+\lx\ly^3+\ly^4}
{(\lx+\ly)(\lx-\ly)^2\lx\ly}\cdot\\
&\qquad\cdot\left(\left(1+3a_1^2\right)\cdot h_{11;\,1}^2
+\left(1+3a_2^2\right)\cdot h_{22;\,2}^2\right)\displaybreak[1]\\
&+2\frac{\lambda_1^2+4\lambda_1\lambda_2+\lambda_2^2}
{(\lambda_1-\lambda_2)^2\lambda_1\lambda_2}\cdot
\left(\lambda_1(1+a_1)^2\cdot h_{11;\,1}^2
+\lambda_2(1+a_2)^2\cdot h_{22;\,2}^2\right)\displaybreak[1]\\
&-\frac{16}{(\lambda_1-\lambda_2)^2(\lambda_1+\lambda_2)}\cdot
\left(\lx(\lx+\ly a_1)(1+a_1)\cdot h_{11;\,1}^2\right.\\
&\quad\left.+\ly(\ly+\lx a_2)(1+a_2)\cdot 
h_{22;\,2}^2\right)\displaybreak[1]\\
=&-4\frac{K^2}H\displaybreak[1]\\
&-2\frac{\left(5\lambda_1^8-4\lambda_1^7\lambda_2
+46\lambda_1^6\lambda_2^2+48\lambda_1^5\lambda_2^3
+72\lambda_1^4\lambda_2^4\right)\lambda_2}
{\left(\lambda_1+\lambda_2\right)\left(\lambda_1-\lambda_2\right)^2
\left(\lambda_1^2+\lambda_1\lambda_2+2\lambda_2^2\right)^2\lambda_1^3}
h_{11;\,1}^2\displaybreak[1]\\
&-2\frac{\left(44\lambda_1^3\lambda_2^5+34\lambda_1^2\lambda_2^6
+8\lambda_1\lambda_2^7+3\lambda_2^8\right)\lambda_2}
{\left(\lambda_1+\lambda_2\right)\left(\lambda_1-\lambda_2\right)^2
\left(\lambda_1^2+\lambda_1\lambda_2+2\lambda_2^2\right)^2\lambda_1^3}
h_{11;\,1}^2\displaybreak[1]\\
&-2\frac{\left(5\lambda_2^8-4\lambda_2^7\lambda_1
+46\lambda_2^6\lambda_1^2+48\lambda_2^5\lambda_1^3
+72\lambda_2^4\lambda_1^4\right)\lambda_1}
{\left(\lambda_2+\lambda_1\right)\left(\lambda_2-\lambda_1\right)^2
\left(\lambda_2^2+\lambda_2\lambda_1+2\lambda_1^2\right)^2\lambda_2^3}
h_{22;\,2}^2\displaybreak[1]\\
&-2\frac{\left(44\lambda_2^3\lambda_1^5+34\lambda_2^2\lambda_1^6
+8\lambda_2\lambda_1^7+3\lambda_1^8\right)\lambda_1}
{\left(\lambda_2+\lambda_1\right)\left(\lambda_2-\lambda_1\right)^2
\left(\lambda_2^2+\lambda_2\lambda_1+2\lambda_1^2\right)^2\lambda_2^3}
h_{22;\,2}^2\displaybreak[1]\\
\le&0.
\end{align*}
We finally apply the maximum principle and our theorem follows.
\end{proof}

\subsection{Direct Consequences}

We obtain a pinching estimate
\begin{corollary}\label{pinching cor}
For a smooth closed strictly convex surface $M_t$ in $\R^3$, 
flowing according to $\dot X=-\A2\nu$, there exists 
$c=c(M_0)$ such that $0<\frac1c\le\frac{\lambda_1}{\lambda_2}\le c$.
\end{corollary}
\begin{proof}
Choose $\epsilon>0$ such that $\lambda_1,\,\lambda_2>\epsilon$
at $t=0$. Theorem \ref{A2 mon thm} and Lemma \ref{A2 conv pre}
imply that 
$$2\epsilon\frac{\left(\frac\lx\ly-1\right)^2}{\frac\lx\ly}
=2\epsilon\frac{(\lx-\ly)^2}{\lx\ly}
\le\frac{(\lx+\ly)(\lx-\ly)^2}{\lx\ly}\le c.$$
We obtain the upper bound on $\frac{\lambda_1}{\lambda_2}$ claimed above. 
Similarly, we obtain an upper bound on $\frac{\lambda_2}{\lambda_1}$.
\end{proof}

Let $\rho_+$ be the minimal radius of enclosing spheres and
$\rho_-$ the maximal radius of enclosed spheres. The quotient
of these radii can be estimated as follows
\begin{corollary}\label{radii cor}
Under the assumptions of Corollary \ref{pinching cor},
$\rho_+/\rho_-$ is bounded above by a constant depending only
on the constant $c(M_0)$ in Corollary \ref{pinching cor}. 
\end{corollary}
\begin{proof}
Combine Corollary \ref{pinching cor}, 
\cite[Theorem 5.1]{AndrewsContractCalc},
and \cite[Lemma 5.4]{AndrewsContractCalc}.
\end{proof}

We also obtain a bound for $\abs{\lx-\ly}$
\begin{corollary}\label{diff bound}
For a smooth closed strictly convex surface $M_t$ in $\R^3$,
flowing according to $\dot X=-\A2\nu$, there exists a
constant $c=c(M_0)$ such that $\abs{\lx-\ly}\le c\cdot\Akl2{1/4}
\le c\cdot\sqrt{H}$.
\end{corollary}
\begin{proof}
This is a direct consequence of Theorem \ref{A2 mon thm}
and Corollary \ref{pinching cor}.
\end{proof}

As in \cite{AndrewsStones}, this estimate on $\abs{\lx-\ly}$ is
``better'' than scaling invariant. It is crucial for the rest
of the proof of Theorem \ref{main thm}.

Let us recall a form of the maximum principle for evolving hypersurfaces.
\begin{lemma}\label{encl lem}
Let $M_t$ and $\tilde M_t$ be two smooth closed strictly convex solutions
to \eqref{flow eqn} on some time interval $\left[0,\,T^*\right)$.
If $M_0$ encloses $\tilde M_0$, then $M_t$ encloses $\tilde M_t$
for any $t\in\left[0,\,T^*\right)$.
\end{lemma}
\begin{proof}
This is a standard consequence of the maximum principle. 
\end{proof}

The next result describes the evolution of spheres. 
\begin{lemma}\label{sphere evol}
Spheres $\partial B_{r(t)}(x_0)$ solve \eqref{flow eqn} for
$t\in[0,\,T)$ with $r(t)=(6(T-t))^{1/3}$ and $T=\frac16r^3(0)$.
\end{lemma}
\begin{proof}
The evolution equation for the radius of a sphere is
$$\dot r(t)=-\frac2{r^2(t)}.$$
\end{proof}

As a consequence, we can estimate the life span of a solution 
in terms of inner and outer radii.
\begin{lemma}\label{life span lem}
Let $\rho_+(t)$ and $\rho_-(t)$ be the inner and outer radii of $M_t$,
respectively. Assume that $M_t$ is a smooth closed strictly convex
solution of \eqref{flow eqn} on a maximal time interval $[0,\,T)$.
Then we have for $t\in[0,\,T)$
$$\frac16\rho_-^3(t)\le T-t\le\frac16\rho_+^3(t).$$
\end{lemma}
\begin{proof}
As $M_t$ contracts to a point, we deduce from Lemma \ref{encl lem}
that $T-t$ is bounded below by the life span of $\partial B_{\rho_-(t)}$
evolving according to \eqref{flow eqn}. So the lower bound follows
from Lemma \ref{sphere evol}. The upper bound is obtained similarly.
\end{proof}

\subsection{Convergence to a Round Point}
We closely follow the corresponding part of \cite{AndrewsStones}.

\begin{proposition}\label{prop int H}
Define $q(t):=\frac1{4\pi}\int\limits_{M_t}KX$.
Then 
$$\left\lvert\langle X-q,\nu\rangle-\frac1{8\pi}
\int\limits_{M_t}H\right\rvert\le\frac1{4\pi}
\cdot\sup_{M_t}\abs{\lx-\ly}\cdot\mathcal H^{2}(M_t),$$
where $\mathcal H^2(M_t)$ denotes the area of $M_t$.
\end{proposition}
\begin{proof}
This is \cite[Proposition 4]{AndrewsStones}.
\end{proof}

We define $r_+(t)$ to be the minimal radius of a sphere, 
centered at $q(t)$, that encloses $M_t$. Similarly, we
define $r_-(t)$ to be the maximal radius of a sphere,
centered at $q(t)$, that is enclosed by $M_t$.

\begin{lemma}\label{radii bound}
Under the assumptions of Theorem \ref{main thm}, 
for $T-t$ sufficiently small, $r_+$ and
$r_-$ are estimated as follows
\begin{align*}
r_+(t)\le&(6(T-t))^{1/3}\cdot\left(1+
c\cdot(T-t)^{1/6}\right),\displaybreak[1]\\
r_-(t)\ge&(6(T-t))^{1/3}\cdot\left(1-
c\cdot(T-t)^{1/6}\right),
\intertext{and}
1\le\frac{r_+}{r_-}\le&1+c\cdot(T-t)^{1/6}.
\end{align*}
\end{lemma}
\begin{proof}
Denote the bounded component of $\R^3\setminus M_t$
by $E_t$. The transformation formula for integrals
implies that
$$\frac1{4\pi}\int\limits_{M_t}KX=\frac1{4\pi}
\int\limits_{\S^2}X\left(\nu^{-1}(\cdot)\right).$$
So we see that $q(t)\in E_t$.
We have 
\begin{align*}
r_+=&\max_{M_t}\,\langle X-q(t),\,\nu\rangle,&
r_-=&\min_{M_t}\,\langle X-q(t),\,\nu\rangle,\\ 
\rho_+=&\min_{p\in\R^3}\max_{M_t}\,\langle X-p,\,\nu\rangle,
\qquad\text{and}&
\rho_-=&\max_{p\in E_t}\min_{M_t}\,\langle X-p,\,\nu\rangle.
\end{align*}

Recall the first variation formula for a vector 
field $Y$ along $M_t$ \cite{TomMCF}
$$\int\limits_{M_t}H
\langle Y,\,\nu\rangle=\int\limits_{M_t}\divergenz_{M_t} Y$$
and get for $p\in E_t$ such that $\rho_+=\max_{M_t}
\langle X-p,\,\nu\rangle$
$$\int\limits_{M_t}H\ge\frac1{\rho_+}\int\limits_{M_t}
H\cdot\langle X-p,\,\nu\rangle=\frac1{\rho_+}\int\limits_{M_t}
\divergence_{M_t}X=\frac1{\rho_+}\int\limits_{M_t}2=\frac2{\rho_+}
\mathcal H^2(M_t).$$
We employ Proposition \ref{prop int H} and deduce that 
\begin{align*}
r_-\ge&\frac1{8\pi}\int\limits_{M_t}H\cdot\left\{1-
2\left(\int_{M_t}H\right)^{-1}\cdot
\sup\limits_{M_t}\abs{\lx-\ly}\cdot\mathcal H^2(M_t)\right\}\\
\ge&\frac1{8\pi}\int\limits_{M_t}H\cdot
\left\{1-\rho_+\cdot\sup\limits_{M_t}\abs{\lx-\ly}\right\}.
\end{align*}
We estimate as follows
\begin{align*}
\rho_+\cdot\sup\limits_{M_t}\abs{\lx-\ly}\le&
c\cdot\rho_+\cdot\Akl2{1/4}&&\text{by Corollary \ref{diff bound}}\\
\le&c\cdot\rho_+\cdot\left(c+\frac c{\rho_-^2}\right)^{1/4}&&
\text{by Lemma \ref{vel bound}}\\
\le&c\cdot(T-t)^{1/6}&&
\end{align*}
by Corollary \ref{radii cor} and Lemma \ref{life span lem}
for $(T-t)$ small. So we obtain
\begin{equation}\label{r- bound}
r_-(t)\ge\frac1{8\pi}\int\limits_{M_t}H\cdot\left(1-
c\cdot(T-t)^{1/6}\right).
\end{equation}
Similar calculations yield
\begin{equation}\label{r+ bound}
r_+(t)\le\frac1{8\pi}\int\limits_{M_t}H\cdot\left(1+
c\cdot(T-t)^{1/6}\right).
\end{equation}
We employ Lemma \ref{life span lem}
$$r_-\le\rho_-\le(6(T-t))^{1/3}\le\rho_+\le r_+$$
and obtain for $(T-t)$ small
$$(6(T-t))^{1/3}\cdot\left(1-c\cdot(T-t)^{1/6}\right)
\le\frac1{8\pi}\int\limits_{M_t}H\le
(6(T-t))^{1/3}\cdot\left(1+c\cdot(T-t)^{1/6}\right).$$
Using \eqref{r- bound} and \eqref{r+ bound} gives the
claimed estimates on $r_-$ and $r_+$, and $r_+/r_-$ is
bounded as stated above.
\end{proof}

\begin{corollary}
Under the assumptions of Theorem \ref{main thm}, we have the
estimate
$$|q(t)-Q|\le c\cdot(T-t)^{1/3+1/6}.$$
Therefore, we obtain the same estimates as in 
Lemma \ref{radii bound}, if we define $r_+$ and $r_-$
using $Q$ instead of $q(t)$. 
\end{corollary}
\begin{proof}
Fix $t_0\in[0,\,T)$. A sphere of radius
$(6(T-t_0))^{1/3}\left(1+c(T-t_0)^{1/6}\right)$, centered
at $q(t_0)$ as defined in Proposition \ref{prop int H},
will enclose $M_t$ for all $t_0\le t<T$. Thus its radius 
at time $t=T$ is an upper bound for $|q(t_0)-Q|$. 
According to Lemma \ref{sphere evol}, the radius of that
sphere evolves as follows
$$r(t)=\left(6\left(\frac16r^3(0)-(t-t_0)\right)\right)^{1/3}.$$
Therefore, we get
\begin{align*}
|q(t_0)-Q|\le&r(T)\umbruch\\
=&\left(6\left(\frac166(T-t_0)\left(1+c(T-t_0)^{1/6}\right)^{3}
-(T-t_0)\right)\right)^{1/3}\umbruch\\
=&(6(T-t_0))^{1/3}\cdot\left(\left(1+c(T-t_0)^{1/6}\right)^3
-1\right)^{1/3}\umbruch\\
\le&(6(T-t_0))^{1/3}\cdot c\cdot(T-t_0)^{1/6}.
\end{align*}
\end{proof}

Next, we want to check, that we can apply a Harnack inequality
\cite[Theorem 5.17]{AndrewsHarnack}. For $F=F(\lambda_i)$,
$\lambda_i>0$, we define
$$\Phi(\kappa_i):=-F\left(\kappa_i^{-1}\right).$$
We say that $\Phi$ is $\alpha$-concave, if $\Phi=\sgn\alpha\cdot
B^\alpha$ for some $B$, where $B$ is positive and concave.
The function $\Phi$ is called the dual function to $F$.

\begin{lemma}\label{alpha concave lem}
The dual function to $F=\A2=\lx^2+\ly^2$ is
$\alpha$-concave for $\alpha\le-2$. 
\end{lemma}
\begin{proof}
We define for $\lambda_i>0$
$$\Phi=\Phi(\lambda_i)=\Phi(h_{ij},\,g_{ij})=
-\left(\frac1{\lx^2}+\frac1{\ly^2}\right)
=\frac{2K-H^2}{K^2}
=\frac{2\det h^i_j-\left(\tr h^i_j\right)^2}{\left(\det h^i_j\right)^2}$$
and obtain
\begin{align*}
\Phi^{ij}\equiv\fracp\Phi{h_{ij}}
=&\frac1{K^2}\left\{\left[2H^2-2K\right]\tilde h^{ij}
-2Hg^{ij}\right\},\displaybreak[1]\\
\Phi^{ij,\,kl}=&\frac1{K^2}\left\{\left[2K-4H^2\right]
\tilde h^{ij}\tilde h^{kl}
+\left[2K-2H^2\right]\tilde h^{ik}\tilde h^{jl}\right.\\
&\quad\left.4H\left[g^{ij}\tilde h^{kl}+g^{kl}\tilde h^{ij}\right]
-2g^{ij}g^{kl}\right\}.
\end{align*}
According to \cite[(5.4)]{AndrewsHarnack}, it suffices to show that
\begin{equation}\label{alpha concave}
\Phi^{ij,\,kl}\eta_{ij}\eta_{kl}\le\frac{\alpha-1}{\alpha\Phi}
\Phi^{ij}\eta_{ij}\Phi^{kl}\eta_{kl}
\end{equation}
for symmetric matrices $(\eta_{ij})$. Terms involving $\eta_{12}$
clearly have the right sign. The remaining terms are a quadratic
form in $(\eta_{11},\,\eta_{22})$. Thus \eqref{alpha concave} is
fulfilled, if
\begin{align*}
\A2\cdot&\left\{\left[2K-4H^2\right]\begin{pmatrix}\frac1{\lx^2}
&\frac1{\lx\ly}\\\frac1{\lx\ly}&\frac1{\ly^2}\end{pmatrix}
+\left[2K-2H^2\right]\begin{pmatrix}\frac1{\lx^2}&0\\
0&\frac1{\ly^2}\end{pmatrix}\right.\\
&\left.+4H\begin{pmatrix}\frac2{\lx}&\frac1{\lx}+\frac1{\ly}\\
\frac1{\lx}+\frac1{\ly}&\frac2{\ly^2}\end{pmatrix}
-2\begin{pmatrix}1&1\\1&1\end{pmatrix}\right\}\le\\
\le&-\frac{\alpha-1}{\alpha}\left\{\left[2H^2-2K\right]\begin{pmatrix}
\frac1{\lx}\\\frac1{\ly}\end{pmatrix}
-2H\begin{pmatrix}1\\1\end{pmatrix}\right\}\otimes
\left\{\left[2H^2-2K\right]\begin{pmatrix}
\frac1{\lx}\\\frac1{\ly}\end{pmatrix}
-2H\begin{pmatrix}1\\1\end{pmatrix}\right\}\\
\end{align*}
or equivalently
\begin{align*}
-6\left(\lx^2+\ly^2\right)\begin{pmatrix}\frac{\ly^2}{\lx^2}&0\\
0&\frac{\lx^2}{\ly^2}\end{pmatrix}\le
-4\frac{\alpha-1}{\alpha}\begin{pmatrix}\frac{\ly^4}{\lx^2}&
\lx\ly\\\lx\ly&\frac{\lx^4}{\ly^2}\end{pmatrix}.
\end{align*}
As $\frac{\alpha-1}{\alpha}\le\frac32$ for $\alpha\le-2$,
we obtain that $\Phi$ is $\alpha$-concave. 
\end{proof}

We are now able to improve our velocity bounds.
\begin{lemma}\label{A2 est}
Under the assumptions of Theorem \ref{main thm}, we obtain 
$$2(6(T-t))^{-2/3}\cdot\left(1-c\cdot(T-t)^{1/12}\right)
\le\A2\le2(6(T-t))^{-2/3}\cdot\left(1+c\cdot(T-t)^{1/12}\right)$$
everywhere on $M_t$ for $(T-t)$ sufficiently small.
\end{lemma}
\begin{proof}
We may assume that $T-t>0$ is so small that we can use the 
results obtained before. Parameterize $M_t$ by $\S^2$ such that
the normal image of $M_t$ at $X(z,\,t)$ equals $z\in\S^2$.
Let us define the support function $s$ of $M_t$ as
$$s(z,\,t):=\langle X(z,\,t),\,z\rangle.$$
Its evolution equation, see e.\,g.\ \cite{AndrewsHarnack}, is
\begin{equation}\label{supp fct evol}
\dt s(z,\,t)=-\A2(z,\,t).
\end{equation}
The $\alpha$-concavity proved in Lemma \ref{alpha concave lem}
allows us to use \cite[Theorem 5.17]{AndrewsHarnack}.
We obtain for $0<t_1<t_2<T$ and $z\in\S^2$, for
two points $(z,\,t_1)$ and $(z,\,t_2)$ with the same normal,
\begin{equation}\label{A2 Harnack}
\frac{\A2(z,\,t_2)}{\A2(z,\,t_1)}\ge
\left(\frac{t_1}{t_2}\right)^{2/3}.
\end{equation}
Let us assume that $q(t)$ is the origin for some fixed
time $t$. As $M_t$ lies between $\partial B_{r_+(t)}(0)$
and $\partial B_{r_-(t)}(0)$, $M_{t+\tau}$ lies outside
$B_{\left(r_-^3(t)-6\tau\right)^{1/3}}(0)$ for any
$0<\tau<T-t$, so
\begin{equation}\label{supp fct bounds}
r_-(t)\le s(\cdot,\,t)\le r_+(t)\quad\text{and}\quad
\left(r_-^3-6\tau\right)^{1/3}\le s(\cdot,\,t+\tau).
\end{equation}
Set $\tau=r_-^{5/2}(t)\cdot(r_+(t)-r_-(t))^{1/2}$ and observe that
$t+\tau<T$, if $(r_+-r_-)^{1/2}\le\frac16r_-^{1/2}$
(by Lemma \ref{life span lem}), or, if
$T-t$ is sufficiently small
(by Lemma \ref{radii bound}). We estimate
\begin{align*}
\A2(z,\,t)\le&\inf\limits_{0\le\tilde\tau\le\tau}
\left\{\left(\frac{t+\tilde\tau}{t}\right)^{2/3}\cdot
\A2(z,\,t+\tilde\tau)\right\}&&\text{by \eqref{A2 Harnack}}
\displaybreak[1]\\
\le&\left(\frac{t+\tau}t\right)^{2/3}\cdot\frac1{\tau}\cdot
\int\limits_t^{t+\tau}\A2(z,\,\tilde\tau)d\tilde\tau&&
\displaybreak[1]\\
\le&\left(1+\frac{\tau}t\right)^{2/3}\cdot\frac1\tau\cdot
(s(z,\,t)-s(z,\,t+\tau))&&\text{by \eqref{supp fct evol}}
\displaybreak[1]\\
\le&\left(1+\frac{\tau}t\right)\cdot\frac1{\tau}\cdot
\left(r_+(t)-\left(r_-^3(t)-6\tau\right)^{1/3}\right)
&&\text{by \eqref{supp fct bounds}}\displaybreak[1]\\
=&\frac{1+\frac\tau t}
{r_-^2\cdot\left(\frac{r_+}{r_-}-1\right)^{1/2}}
\cdot\left(\frac{r_+}{r_-}-\left(1
-6\cdot\left(\frac{r_+}{r_-}-1\right)^{1/2}\right)^{1/3}\right).&&
\end{align*}
The maximal existence time $T$ is bounded below in terms of the
radius of a sphere enclosed by $M_0$. So we may also assume
that $t$ is bounded below by a positive constant. A very crude
estimate gives 
$$\tau\le r_-^{5/2}\cdot r_+^{1/2}
\le r_+^3\le c\cdot(T-t),$$ 
so we obtain
$$\frac{1+\frac\tau t}{r_-^2}\le
(6(T-t))^{-2/3}\cdot\left(1+c\cdot(T-t)^{1/6}\right).$$
Observe that for $\abs x\le\frac12$, we have
$$-(1-x)^{1/3}\le-1+\tfrac13x+\tfrac13x^2\quad
\text{and }\quad
(1+x)^{1/3}\ge1+\tfrac13x-\tfrac13x^2.$$
We conclude for small $(T-t)$
\begin{align*}
\A2(z,\,t)\le&(6(T-t))^{-2/3}\cdot\left(1+c\cdot(T-t)^{1/6}\right)
\cdot\left(2+13\cdot\left(\frac{r_+}{r_-}-1\right)^{1/2}\right)\\
\le&2\cdot(6(T-t))^{-2/3}\cdot\left(1+c\cdot(T-t)^{1/12}\right).
\end{align*}
For the lower bound on $\A2$, we proceed similarly and use
$\tau=r_-^{5/2}(t)\cdot(r_+(t)-r_-(t))^{1/2}$, 
$$\left(r_-^3+6\tau\right)^{1/3}\le s(z,\,t-\tau)\quad\text{and}
\quad s(z,\,t)\le r_+(t).$$
\end{proof}

We have the following estimate for the principal curvatures
\begin{lemma}\label{lambda bounds}
Under the assumptions of Theorem \ref{main thm}, we obtain
$$(6(T-t))^{-1/3}\cdot\left(1-c\cdot(T-t)^{1/12}\right)\le
\lx,\,\ly\le(6(T-t))^{-1/3}\cdot\left(1+c\cdot(T-t)^{1/12}\right)$$
on $M_t$ for small $(T-t)$. 
\end{lemma}
\begin{proof}
As $H^2=2\A2-(\lx-\ly)^2$, we obtain
\begin{align}\label{lx expr}
\begin{split}
\lx=&\tfrac12(\lx+\ly)+\tfrac12(\lx-\ly)\\
=&\tfrac12\sqrt{2\A2-(\lx-\ly)^2}+\tfrac12(\lx-\ly).
\end{split}
\end{align}
Combining Lemmata \ref{diff bound} and \ref{A2 est}, 
we get $\abs{\lx-\ly}\le c\cdot(T-t)^{-1/6}$. We
use Lemma \ref {A2 est} and \eqref{lx expr}. The 
claimed inequality follows.
\end{proof}

\begin{proof}[Proof of Theorem \ref{main thm}:]
Lemma \ref{lambda bounds} implies, that, everywhere on $M_t$, 
the quotient $\lx/\ly$ tends to $1$ as $t\uparrow T$. Then 
we can apply known results, see e.\,g.\ 
\cite[Theorem 2]{Andrews2dnonconcave}, to conclude that 
the rescaled surfaces converge smoothly to the unit sphere
$\S^2\subset\R^3$.
\end{proof}

A standard way of rescaling \cite{AndrewsContractCalc} 
is to consider the embeddings $\tilde X(\cdot,\,t)$,
$$\tilde X(z,\,t):=(6(T-t))^{-1/3}(X(z,\,t)-Q)$$ 
with $Q$ as in Theorem \ref{main thm}. Define 
the time function $\tau(t):=\frac16\log T-\frac16\log(T-t)$.
Then we have, using suggestive notation, 
the following evolution equation
$$\frac d{d\tau}\tilde X=-\abs{\tilde A}^2\tilde\nu+2\tilde X$$
and our a priori estimates imply, that, for $\tau\to\infty$, 
$\tilde M_t$ converges exponentially to $\S^2$.


\section{Finding Monotone Quantities}\label{algor sec}

\subsection{The Algorithm}
We use a sieve algorithm and start with symmetric rational functions 
of the principal curvatures
as candidates for test functions, e.\,g.\ 
$$w=\frac{p_1(\lx,\,\ly)}{p_2(\lx,\,\ly)}
=\frac{(\lx+\ly)(\lx-\ly)^2}{\lx\ly}.$$
Here, $p_1\neq0$ and $p_2\neq0$ are homogeneous polynomials. 

In the end, we want to find functions $w$ such that 
$W:=\sup_{M_t}w$ is monotone and ensures convergence to round
spheres. 

We check, whether these test functions $w$ fulfill the
following conditions. 
\begin{enumerate}
\item\label{step one} 
 \begin{enumerate}
 \item $p_1(\lx,\,\ly),\,p_2(\lx,\,\ly)\ge0$ for $0<\lx,\,\ly$,
 \item $p_1(\lx,\,\ly)=0$ for $\lx=\ly>0$.
 \end{enumerate}
\item\label{step two} $\deg p_1>\deg p_2$.
\item\label{step three} $\fracp{w(1,\ly)}{\ly}<0$ for $0<\ly<1$ and 
  $\fracp{w(1,\ly)}{\ly}>0$ for $\ly>1$.
\item\label{step four} $\dt w-F^{ij}w_{;\,ij}\le0$ 
 \begin{enumerate}
 \item\label{step four a} for terms without derivatives of $(h_{ij})$,
 \item\label{step four b} for terms involving derivatives of $(h_{ij})$,
   if $w_{;\,i}=0$ for $i=1,\,2$.
 \end{enumerate}
\end{enumerate}

\subsection{Motivation and Randomized Tests}
We restrict 
our attention to non-negative polynomials $p_i$. For all flow equations 
considered, spheres contract to points and stay spherical. So
we can only find monotone quantities, if $\deg p_1\le\deg p_2$
or $p_1(\lambda,\,\lambda)=0$. 

If $\deg p_1<\deg p_2$,
we obtain that $W$ is decreasing on any self-similarly
shrinking surface. So this does not imply convergence to a
sphere. The counterexamples in \cite{Andrews2dnonconcave} show 
for normal velocities of homogeneity larger than $1$, that 
the pinching ratio $\sup_{M_t}\ly/\lx$ (for $\ly>\lx$) 
will increase 
during the flow for appropriate initial surfaces. Therefore,
we require in step \eqref{step two}, that $\deg p_1>\deg p_2$. 

Condition \eqref{step three} ensures that the quantity decreases, if 
the eigenvalues approach each other. This excludes 
especially local zeroes of $w(1,\,\ly)$ for $\ly\neq 1$. 

In all these steps, 
inequalities are tested by evaluating both sides at
random numbers. If an inequality is not violated for 
sufficiently many tuples of random numbers, we move to the
next step and keep the candidate, otherwise we start with
another candidate. 

The hard part is to test, whether $\dt w-F^{ij}w_{;\,ij}\le0$
holds at a point, where $w_{;\,i}=0$. 

We assume in step \eqref{step four a}, 
that $h_{ij;\,k}=0$. 
Recall the algebraic fact, that all products of $H$ and $\A2$
of given homogeneity form a basis for the symmetric homogeneous
polynomials of that degree. We represent the polynomials
$p_i$ in this basis. At random values for
$\lx,\,\ly$, we compute $\oprokl H$ and $\opr{\A2}$.
This can be combined according to the rules of 
differentiation and yields $\oprokl w$,
evaluated at these numbers.

If not all components of $h_{ij;\,k}$ vanish (step 
\eqref{step four b}), we also have to
choose random numbers for $h_{ij;\,k}$ that fulfill the
extremal condition $w_{;\,i}=0$. As above, we evaluate 
$\oprokl w$ at the random numbers chosen.
Here we can ignore all terms that do not contain derivatives 
of the second fundamental form. The evaluation of the
remaining terms is more involved than in the last step,
but can be done similarly according to the various rules
of differentiation. 

Now we iterate steps \eqref{step four a} and 
\eqref{step four b}. If all tests yield 
$\oprokl w\le0$ at critical points of $w$, it seems likely 
that we have found an appropriate test quantity. Indeed,
if this inequality is fulfilled, the maximum principle 
implies that $W=\max_{M_t}w$ is non-increasing in time.

We implemented this algorithm in a C-program and used it
to find all the new test functions of this paper.
 
Obviously, the computing time depends on the
number of tests performed. The following computing times
are measured for the quantity 
\eqref{A2 quantity} on an Intel Pentium 4, 2.4 GHz, running
Linux 2.4.24 and GNU C-compiler 2.95.4. 
The number of tests per second
for step \eqref{step four a} is $1.6\cdot 10^5$ 
and $5.8\cdot 10^3$ for step \eqref{step four b}.
The other steps are comparable to step \eqref{step four a}.
In steps \eqref{step one} to \eqref{step three}, 
the calculations do not depend on the normal velocity.

It seems worth noting, that, after testing with enough random
numbers in an appropriate range, every candidate for a monotone
quantity that we checked, turned out to be a useful test quantity.
In that sense, algorithm and program seem to be correct.

We are convinced that it is possible to implement this
algorithm for surfaces without using random numbers. Evolution
equations can be computed algebraically and Sturm's algorithm
can be used to test for non-negativity.

We expect that similar algorithms will be used to find (monotone)
test functions for other (geometric) problems. 

\section{$H^3$-Flow}

In this and the following sections, we will consider strictly convex
surfaces contracting according to $\dt X=-F\nu$ for several normal
velocities $F$. We will not repeat parts of the argument that are
very similar to the respective parts in the proof for $F=\A2$.

As the theorems for these flow equations agree essentially with
Theorem \ref{main thm}, we will state them in concise form as
follows.

\begin{theorem}\label{H3 thm}
A smooth closed strictly convex surface in $\R^3$, contracting
with normal velocity $H^3$, converges to a round point in finite time.
\end{theorem}

\subsection{A Monotone Quantity}

\begin{theorem}
For a family of smooth closed strictly convex surfaces $M_t$ 
in $\R^3$, flowing according to $\dt X=-H^3\nu$, 
$$\max\limits_{M_t}
\frac{\left(\lx^2+\lx\ly+\ly^2\right)(\lx+\ly)^2(\lx-\ly)^2}
{2\left(\lx^2-\lx\ly+\ly^2\right)\lx\ly}$$
is non-increasing in time.
\end{theorem}
\begin{proof}
We compute
\begin{align*}
F=&H^3,\displaybreak[1]\\
F^{ij}=&3H^2g^{ij},\displaybreak[1]\\
\oprokl H=&H^3\A2+6H\sum\hii,\displaybreak[1]\\
\opr{\A2}=&6H^2\Akl22-4H^3\A3\\
&-6H^2\sum\hij+12H\sum\lk\hii,\displaybreak[1]\\
\tilde w=&\log\left(\frac{\left(\lx^2+\lx\ly+\ly^2\right)
(\lx+\ly)^2(\lx-\ly)^2}
{2\left(\lx^2-\lx\ly+\ly^2\right)\lx\ly}\right)\\
=&\log\left(\frac{-H^6+\A2H^4+2\Akl22H^2}{-H^4+4\A2H^2-3\Akl22}\right).
\end{align*}
In a critical point of $\tilde w$, we obtain
\begin{align*}
h_{22;\,1}=&\frac{\ly}{\lx}\cdot
\frac{3\lx^6-2\lx^5\ly+4\lx^4\ly^2+4\lx^3\ly^3+4\lx^2\ly^4-2\lx\ly^5+\ly^6}
{3\ly^6-2\ly^5\lx+4\ly^4\lx^2+4\ly^3\lx^3+4\ly^2\lx^4-2\ly\lx^5+\lx^6}
\cdot h_{11;\,1},
\end{align*}
\begin{align*}
\oprokl{\tilde w}=&\left(\tfrac{-6H^5+4\A2H^3+4\Akl22H}
{-H^6+\A2H^4+2\Akl22H^2}-\tfrac{-4H^3+8\A2H}{-H^4+4\A2H^2-3\Akl22}
\right)\cdot\\
&\quad\cdot\left(\A2H^3+6H\sum\hii\right)\displaybreak[1]\\
&+\left(\tfrac{H^4+4\A2H^2}{-H^6+\A2H^4+2\Akl22H^2}
-\tfrac{4H^2-6\A2}{-H^4+4\A2H^2-3\Akl22}\right)\cdot\\
&\quad\cdot\left(6\Akl22H^2-4H^3\A3\right.\\
&\qquad\left.-6H^2\sum\hij+12H\sum\lk\hii\right)\displaybreak[1]\\
&+\left(\tfrac{-30H^4+12\A2H^2+4\Akl22}{-H^6+\A2H^4+2\Akl22H^2}
-\tfrac{-12H^2+8\A2}{-H^4+4\A2H^2-3\Akl22}\right)\cdot\\
&\quad\cdot\left(-3H^2\sum\hii\right)\displaybreak[1]\\
&+\left(\tfrac{4H^2}{-H^6+\A2H^4+2\Akl22H^2}
-\tfrac{-6}{-H^4+4\A2H^2-3\Akl22}\right)\cdot\\
&\quad\cdot\left(-12H^2\sum\li\lj\hii\right)\displaybreak[1]\\
&+\left(\tfrac{4H^3+8H\A2}{-H^6+\A2H^4+2\Akl22H^2}
-\tfrac{8H}{-H^4+4\A2H^2-3\Akl22}\right)\cdot\\
&\quad\cdot\left(-6H^2\sum(\li+\lj)\hii\right)\displaybreak[1]\\
=&-12\frac{(\lx+\ly)^2\lx^3\ly^3}
{\left(\lx^2-\lx\ly+\ly^2\right)
\left(\lx^2+\lx\ly+\ly^2\right)}\displaybreak[1]\\
&-\frac{12(\lx+\ly)^4\ly}{\left(\lx^6-2\lx^5\ly+4\lx^4\ly^2
+4\lx^3\ly^3+4\lx^2\ly^4-2\lx\ly^5+3\ly^6\right)^2}\cdot\displaybreak[1]\\
&\quad\cdot\frac1{\left(\lx^2-\lx\ly+\ly^2\right)
\left(\lx^2+\lx\ly+\ly^2\right)(\lx-\ly)^2\lx^3}\cdot\displaybreak[1]\\
&\quad\cdot\left(2\lx^{16}-10\lx^{15}\ly+10\lx^{14}\ly^2
+42\lx^{13}\ly^3-135\lx^{12}\ly^4\right.\\
&\qquad\left.+274\lx^{11}\ly^5-272\lx^{10}\ly^6+258\lx^9\ly^7
-144\lx^8\ly^8+262\lx^7\ly^9\right.\\
&\qquad\left.-166\lx^6\ly^{10}+98\lx^5\ly^{11}+6\lx^4\ly^{12}
-20\lx^3\ly^{13}+14\lx^2\ly^{14}\right.\\
&\qquad\left.-4\lx\ly^{15}+\ly^{16}\right)
\cdot h_{11;\,1}^2\displaybreak[1]\\
&-(\ldots)\cdot h_{22;\,2}^2\displaybreak[1]\\
\le&0.
\end{align*}
We applied Sturm's algorithm \cite{SturmThm} 
to obtain the last inequality.

Here and in the rest of the paper, we have sometimes used
a computer algebra program for the calculations involving
longer polynomials. Moreover, we use $(\ldots)$ in 
$(\ldots)\cdot h_{22;\,2}^2$ to abbreviate terms that 
are, up to interchanging $\lx$ and $\ly$,
equal to the respective factors in front of $h_{11;\,1}^2$. 

We have applied the following two identities in order to rewrite
terms involving derivatives of the second fundamental form
\begin{align*}
\sum\lambda_i^\alpha\lambda_j^\beta\lambda_k^\gamma\hij=&
\left(\lx^\alpha\lx^\beta\lx^\gamma
+\lx^\alpha\ly^\beta\ly^\gamma a_1^2
+\ly^\alpha\lx^\beta\ly^\gamma a_1^2
+\ly^\alpha\ly^\beta\lx^\gamma a_1^2\right)\cdot h_{11;\,1}^2\\
&+(\ldots)\cdot h_{22;\,2}^2,\displaybreak[1]\\
\sum\lambda_i^\alpha\lambda_j^\beta\lambda_k^\gamma\hii=&
\left(\lx^\alpha+\ly^\alpha a_1\right)
\left(\lx^\beta+\ly^\beta a_1\right)\lx^\gamma h_{11;\,1}^2
+(\ldots)\cdot h_{22;\,2}^2.
\end{align*}
They hold for all $\alpha,\,\beta,\,\gamma\in\R$.
\end{proof}

\subsection{Velocity Bounds}

The following lemma is known, see \cite{FelixMZ}. Nevertheless, 
we include it, as we will use some of the 
calculations of its proof later on.
\begin{lemma}\label{H lower bound}
For a family of smooth closed strictly convex surfaces 
$M_t\subset\R^3$, $0\le t<T$, flowing according to
$\dt X=-F\nu$ with $F=H^\alpha$, $\alpha>1$, a positive
lower bound on the principal curvatures, $\lx,\,\ly\ge\epsilon
>0$, is preserved during the evolution.
\end{lemma}
\begin{proof}
Combining \eqref{h evol}, \eqref{Riem}, and \eqref{Ricci}
yields
\begin{align*}
\dt h_{ij}-F^{kl}h_{ij;\,kl}=&F^{kl}h^a_kh_{al}\cdot h_{ij}
-F^{kl}h_{kl}\cdot h^a_ih_{aj}\\
&-Fh^k_ih_{kj}+F^{kl,\,rs}h_{kl;\,i}h_{rs;\,j}.
\end{align*}
We wish to apply the maximum principle for tensors
\cite{ChowLuMP,HamiltonFour86,HamiltonThree}. So we define
$$M_{ij}=h_{ij}-\epsilon g_{ij},$$
use \eqref{g evol} and obtain
\begin{align*}
\dt M_{ij}-F^{kl}M_{ij;\,kl}=&F^{kl}h^a_kh_{al}\cdot h_{ij}
-F^{kl}h_{kl}h^a_ih_{aj}-Fh^k_ih_{kj}\\
&+F^{kl,\,rs}h_{kl;\,i}h_{rs;\,j}+2\epsilon Fh_{ij}.
\end{align*}
We now specialize to the normal velocity $F=H^\alpha$. 
It is easy to see that 
$$F^{kl,\,rs} h_{kl;\,i} h_{rs;\,j}\ge0.$$
We have to test the zero eigenvalue condition. 
Assume that $\xi$ is a zero eigenvalue of $(M_{ij})$,
$h_{ij}\xi^j=\epsilon g_{ij}\xi^j$, with $g_{ij}\xi^i\xi^j=1$.
We may assume, that in our coordinate system, we have $\xi=(1,\,0)$
and $(h_{ij})=\left(\begin{smallmatrix}\epsilon&0\\0&\lambda
\end{smallmatrix}\right)$ with $\lambda\ge\epsilon>0$.
Direct calculations yield 
\begin{align*}
\left(\dt M_{ij}-F^{kl}M_{ij;\,kl}\right)\xi^i\xi^j\ge&
\alpha H^{\alpha-1}\A2\epsilon-\alpha H^\alpha\epsilon^2
-H^\alpha\epsilon^2+2\epsilon H^\alpha\epsilon
\displaybreak[1]\\
=&H^{\alpha-1}\left(\alpha\epsilon\lambda(\lambda-\epsilon)
+\epsilon^3+\epsilon^2\lambda\right)>0.
\end{align*}
The lemma follows from the maximum principle.
\end{proof}

The next lemma appears also in \cite{FelixMZ}. Once again,
we will use the following calculations later on.
\begin{lemma}\label{H upper bound}
For a family of closed smooth strictly convex surfaces 
$M_t\subset\R^3$, $0\le t<T$, flowing according to
$\dt X=-F\nu$ with $F=H^\beta$, $\beta>0$, the velocity
$F$ is bounded as in Lemma \ref{vel bound},
in terms of $\beta$, the initial data, and the
radius $R$ of an enclosed sphere. 
\end{lemma}
\begin{proof}
We proceed as in Lemma \ref{vel bound}, use $\alpha=\frac12R$ and
obtain, see e.\,g. \cite[Lemma 5.4]{OSMathZ}, \cite{OSKnutCrelle},
\begin{align*}
\fracp{F}{g_{kl}}=&-F^{il}h^k_i,\displaybreak[1]\\
\oprokl F=&FF^{ij}h^k_ih_{kj},\displaybreak[1]\\
\dt X^\alpha-F^{ij}X^\alpha_{;\,ij}=&\left(-F+F^{ij}h_{ij}\right)
\nu^\alpha,\displaybreak[1]\\
\dt\nu^\alpha-F^{ij}\nu^\alpha_{;\,ij}=&F^{ij}h^k_ih_{kj}\nu^\alpha,
\displaybreak[1]\\
\oprokl{\langle X,\,\nu\rangle}=&
-F-F^{ij}h_{ij}+F^{ij}h^k_ih_{kj}\langle X,\,\nu\rangle,
\displaybreak[1]\\
\dt\log\left(\frac F{\langle X,\,\nu\rangle-\alpha}\right)\qquad&\\
-F^{ij}\left(\log\left(\frac F{\langle X,\,\nu\rangle-\alpha}\right)
\right)_{;\,ij}=&\frac1{\langle X,\,\nu\rangle-\alpha}
\left(F+F^{ij}h_{ij}-\alpha F^{ij}h^k_ih_{kj}\right).
\end{align*}
So far, we did not use the fact, that $F=H^\beta$.
In an increasing maximum of 
$\frac{F}{\langle X,\,\nu\rangle-\alpha}$, we get the inequality
$F+F^{ij}h_{ij}-\alpha F^{ij}h^k_ih_{kj}\ge0$ and deduce there 
$$\frac1\alpha\frac{1+\beta}\beta\ge\frac{\A2}H\ge\tfrac12H.$$
Our Lemma follows.
\end{proof}

\subsection{Concavity}

\begin{lemma}
The dual function to 
$F=H^\alpha$, $\alpha>0$, defined before
Lemma \ref{alpha concave lem}, is $-\alpha$-concave.
\end{lemma}
\begin{proof}
We use the same notation as before and obtain
\begin{align*}
\Phi=&-H^\alpha K^{-\alpha},\displaybreak[1]\\
\Phi^{ij}=&-\alpha H^{\alpha-1}K^{-\alpha}g^{ij}
+\alpha H^\alpha K^{-\alpha}\tilde h^{ij},\displaybreak[1]\\
\Phi^{ij,\,kl}=&-\alpha(\alpha-1)H^{\alpha-2}K^{-\alpha}g^{ij}g^{kl}\\
&+\alpha^2H^{\alpha-1}K^{-\alpha}g^{ij}\tilde h^{kl}
+\alpha^2H^{\alpha-1}K^{-\alpha}\tilde h^{ij}g^{kl}\\
&-\alpha^2H^\alpha K^{-\alpha}\tilde h^{ij}\tilde h^{kl}
-\alpha H^\alpha K^{-\alpha}\tilde h^{ik}\tilde h^{jl}.
\end{align*}
We have to prove that
$$\Phi^{ij,\,kl}\eta_{ij}\eta_{kl}\le
\frac{\alpha+1}{\alpha\Phi}\Phi^{ij}\eta_{ij}\Phi^{kl}\eta_{kl}$$
or 
\begin{align*}
0\le&-(\alpha+1)\begin{pmatrix}\frac{\ly^2}{\lx^2}&1\\
1&\frac{\lx^2}{\ly^2}\end{pmatrix}
+(\alpha-1)\begin{pmatrix}1&1\\1&1\end{pmatrix}\\
&-\alpha H\begin{pmatrix}\frac2\lx&\frac1\lx+\frac1\ly\\
\frac1\lx+\frac1\ly&\frac2\ly\end{pmatrix}
+\alpha H^2\begin{pmatrix}\frac1{\lx^2}&\frac1{\lx\ly}\\
\frac1{\lx\ly}&\frac1{\ly^2}\end{pmatrix}\\
&+H^2\begin{pmatrix}\frac1{\lx^2}&0\\0&\frac1{\ly^2}\end{pmatrix}
\displaybreak[1]\\
=&\begin{pmatrix}2\frac\ly\lx&-2\\-2&2\frac\lx\ly\end{pmatrix}.
\end{align*}
\end{proof}

\subsection{Some Constants}
We wish to obtain precise bounds on the normal
velocity $F$ near $t=T$. This proof is almost identical for all our
test functions in Table \ref{table} with a factor $\lx\ly$
in the denominator. So it seems appropriate to state this proof
only once with appropriate constants depending on the normal
velocity and the test function. These constants are
\begin{itemize}
\item $c_h$: the homogeneity of $F$ in terms of the principal
  curvatures.
\item $c_1$: the value of $F$ at $\lx=\ly=1$.
\item $c_\alpha$: a positive constant, such that the dual function
  to $F$, defined before Lemma \ref{alpha concave lem},
  is a $-c_\alpha$-concave function. It turns out, that, for
  all flows considered here, we can choose $c_\alpha=c_h$.
\item $c_d$: a constant depending on the difference of the 
  degrees of the numerator $d_n$ and the denominator $d_d$
  of the test function $w$,
  $c_d:=\frac12(2-d_n+d_d)$. This constant is defined such that
  \begin{equation}\label{diff est}
  \abs{\lx-\ly}\le c\cdot H^{c_d}
  \end{equation} for a pinched surface for which $\max_{M_t}w$
  is non-increasing. 
\end{itemize}

For the flow equations considered here, these constants are as
in Table \ref{c tab}. We assume there that $\alpha\ge2$ and
$\beta\ge0$.

It is important for us that $c_d<1$ as it implies that 
Inequality \eqref{diff est} is not scaling invariant.
\begin{table}
\def\platzz{\raisebox{0em}[1.2em][0.7em]{\rule{0em}{2em}}}
$$\begin{array}{|c||c|c|c|c|c|c|c|}\hline
F & \platzz\B{\alpha} & H^2 & H^3 & H^4 & H\A2 &
\abs{A}^4 & \A2+\beta H^2\\\hline\hline
c_h & \platzz\alpha & 2 & 3 & 4 & 3 & 4 & 2\\\hline
c_1 & \platzz 2 & 4 & 8 & 16 & 4 & 4 & 2+4\beta\\\hline
c_\alpha & \platzz\alpha & 2 & 3 & 4 & 3 & 4 & 2\\\hline
c_d & \platzz\tfrac12(3-\alpha) & \tfrac12 & 0 & -2 & 0 
& -\tfrac12 & \tfrac12\\\hline
\frac{1-c_d}{1+c_h} & \platzz \frac{\alpha-1}{2(\alpha+1)}
&\frac16 &\frac14 &\frac35 &\frac14 & \tfrac3{10} &\frac16\\\hline
\end{array}$$
\caption{Some constants}\label{c tab}
\end{table}

\subsection{Pinching}
We show that our surfaces are pinched during the evolution,
i.\,e.\ that there exists a constant $c>0$, depending on our test
quantity, especially on the upper bound for it, and on the positive
lower bound for the principal curvatures, $\epsilon$, such that
$$0<\frac1c\le\frac\li\lj\le c\quad\text{for }i,\,j\in\{1,\,2\}.$$

The following proof does not apply directly to the case $F=K$ 
considered in \cite{AndrewsStones}, but the result is also true
in that case.

By direct inspection, we see that all our test quantities $w$ are
such that $w\cdot\frac{\lx\ly}{(\lx-\ly)^2}$ is bounded below by
a positive constant, depending especially on $\epsilon$. So we see
that $$\frac{\left(\frac\lx\ly-1\right)^2}{\frac\lx\ly}=
\frac{(\lx-\ly)^2}{\lx\ly}\le c.$$
Thus $\frac\lx\ly$ is bounded above and the surface is pinched.

\subsection{Evolution of Spheres}
The radius of contracting spheres fulfills the ordinary 
differential equation
$$\dot r(t)=-\frac{c_1}{r(t)^{c_h}}.$$
A solution is given by
$$r(t)=(c_1\cdot(1+c_h)\cdot(T-t))^{\fracm1{1+c_h}},$$
so that inner and outer radii are related to the life span $T-t$
of $M_t$ as follows
$$\rho_-\le(c_1\cdot(1+c_h)\cdot(T-t))^{\fracm1{1+c_h}}\le
\rho_+(t).$$

\subsection{Bounds for Radii}
In order to prove bounds for the radii $_+$ and $r_-$, we
proceed as in the proof of Lemma \ref{radii bound} and use
$\abs{\lx-\ly}\le c\cdot H^{c_d}$
$$r_-\ge\frac1{8\pi}\int\limits_{M_t}H\cdot\left\{
1-c\cdot\rho_+\cdot H^{c_d}\right\}.$$
For each $F$ in Table \ref{c tab}, there exists $c_F>0$
such that
$$0<\frac1{c_F}H\le F^{1/c_h}\le c_FH,$$
so we can apply the variants of Lemma \ref{vel bound} for other 
curvature functions. For $t$ close to $T$, we get
\begin{align*}
r_-\ge&\frac1{8\pi}\int\limits_{M_t}H\cdot
\left\{1-c\cdot(T-t)^{\fracm{1-c_d}{1+c_h}}\right\}\displaybreak[1]
\intertext{and similarly, we obtain}
r_+\le&\frac1{8\pi}\int\limits_{M_t}H\cdot
\left\{1+c\cdot(T-t)^{\fracm{1-c_d}{1+c_h}}\right\}.
\end{align*}
Then we get for small $T-t$
\begin{align*}
&(c_1\cdot(1+c_n)\cdot(T-t))^{\fracm1{1+c_h}}\cdot
\left\{1-c\cdot(T-t)^{\fracm{1-c_d}{1+c_h}}\right\}\le\\
\le\frac1{8\pi}\int\limits_{M_t}H\le&
(c_1\cdot(1+c_h)\cdot(T-t))^{\fracm1{1+c_h}}\cdot
\left\{1+c\cdot(T-t)^{\fracm{1-c_d}{1+c_h}}\right\},\displaybreak[1]\\
r_-\ge&(c_1\cdot(1+c_h)\cdot(T-t))^{\fracm1{1+c_h}}\cdot
\left\{1-c\cdot(T-t)^{\fracm{1-c_d}{1+c_h}}\right\},\displaybreak[1]\\
r_+\le&(c_1\cdot(1+c_h)\cdot(T-t))^{\fracm1{1+c_h}}\cdot
\left\{1+c\cdot(T-t)^{\fracm{1-c_d}{1+c_h}}\right\},\displaybreak[1]
\intertext{and}
1\le&\frac{r_+}{r_-}\le1+c\cdot(T-t)^{\fracm{1-c_d}{1+c_h}}.
\end{align*}

\subsection{Precise Velocity Bounds}
We use the notation of Lemma \ref{A2 est}.
The Harnack inequality \cite[Theorem~5.17]{AndrewsHarnack}
implies for $t_2>t_1>0$
$$\frac{F(z,\,t_2)}{F(z,\,t_1)}\ge
\left(\frac{t_1}{t_2}\right)^{\fracm{c_\alpha}{1+c_\alpha}}.$$
Note that spheres evolve such that we get
$$\left(r_-(t)^{1+c_h}-c_1(1+c_h)\tau\right)^{\fracm1{1+c_h}}
\le s(\cdot,\,t+\tau).$$
We set 
$$\tau:=r_-(t)^{1+c_h}\cdot\left(\frac{r_+(t)}{r_-(t)}-1\right)^{1/2}$$
and get
\begin{align*}
F(z,\,t)\le&\frac{1+c\cdot(T-t)}{r_-(t)^{c_h}
\left(\frac{r_+(t)}{r_-(t)}-1\right)^{1/2}}\cdot\\
&\cdot\left(\frac{r_+(t)}{r_-(t)}-\left(1-c_1\cdot(1+c_h)\cdot
\left(\frac{r_+(t)}{r_-(t)}-1\right)^{1/2}\right)^{\fracm1{1+c_h}}\right).
\end{align*}
Use for $0\le x\le c(c_h)$
$$-(1-x)^{\fracm1{1+c_h}}\le-1+\frac1{1+c_h}x+\frac1{1+c_h}x^2.$$
We get
\begin{align}\label{F bound}
F(z,\,t)\le&\frac{1+c\cdot(T-t)}{r_-(t)^{c_h}}
\cdot\left(c_1+c\left(\frac{r_+}{r_-}-1\right)^{1/2}\right)
\displaybreak[1]\\
\le&c_1\cdot\left(c_1\cdot(1+c_h)\cdot(T-t)\right)^{-\frac{c_h}{1+c_h}}
\cdot\left(1+c\cdot(T-t)^{\fracm12\cdot\fracm{1-c_d}{1+c_h}}\right)
\nonumber
\end{align}
and a similar lower bound follows. 

\subsection{Convergence of Principal Curvatures}
We consider $F=H^\alpha$ and obtain
\begin{align*}
\lx=&\tfrac12(\lx+\ly)+\tfrac12(\lx-\ly)\displaybreak[1]\\
=&\tfrac12 F^{1/\alpha}+\tfrac12(\lx-\ly)\displaybreak[1]\\
\le&\left(2^\alpha\cdot(1+\alpha)
\cdot(T-t)\right)^{-\fracm1{1+\alpha}}\cdot
\left(1+c\cdot(T-t)^{\fracm12\cdot\fracm{1-c_d}{1+\alpha}}\right).
\end{align*}
A similar lower bound is proved analogously.
Theorem \ref{H3 thm} follows.


\section{$H^2$-Flow}

\begin{theorem}\label{H2 thm}
A smooth closed strictly convex surface in $\R^3$, contracting
with normal velocity $H^2$, converges to a round point in finite time.
\end{theorem}

\begin{theorem}
For a family of smooth closed strictly convex surfaces $M_t$ 
in $\R^3$, flowing according to $\dt X=-H^2\nu$, 
$$\max\limits_{M_t}\frac{(\lx+\ly)^3(\lx-\ly)^2}
{2\left(\lx^2+\ly^2\right)\lx\ly}$$
is non-increasing in time.
\end{theorem}
\begin{proof}
We set 
\begin{align*}
w=&\frac{(\lx+\ly)^3(\lx-\ly)^2}
{2\left(\lx^2+\ly^2\right)\lx\ly}
=\frac{H^3\left(2\A2-H^2\right)}
{\A2\left(H^2-\A2\right)}.
\end{align*}
In a critical point of $\tilde w$, we get, based on
computer algebra calculations, 
\begin{align*}
\oprokl w=&-2\frac{(\lx+\ly)^4(\lx-\ly)^2\lx\ly}
{\left(\lx^2+\ly^2\right)^2}\umbruch\\
&-\frac{(\lx+\ly)^4}{\left(\lx^2+\ly^2\right)^2
\left(\lx^4-\lx^3\ly+7\lx^2\ly^2-\lx\ly^3+2\ly^4\right)^2
\lx^4}\cdot\umbruch\\
&\quad\cdot\left(5\lx^{12}-24\lx^{11}\ly+112\lx^{10}\ly^2
-164\lx^9\ly^3+529\lx^8\ly^4\right.\\
&\qquad\left.-448\lx^7\ly^5+952\lx^6\ly^6-312\lx^5\ly^7
+391\lx^4\ly^8-72\lx^3\ly^9\right.\\
&\qquad\left.+56\lx^2\ly^{10}-4\lx\ly^{11}
+3\ly^{12}\right)\cdot h_{11;\,1}^2
+(\ldots)\cdot h_{22;\,2}^2\umbruch\\
\le&0.
\end{align*}
We apply the maximum principle. Our claim follows.
\end{proof}
Theorem \ref{H2 thm} follows. 


\section{$H^4$-Flow}

\begin{theorem}\label{H4 thm}
A smooth closed strictly convex surface in $\R^3$, contracting
with normal velocity $H^4$, converges to a round point in finite time.
\end{theorem}

\begin{theorem}
For a family of smooth closed strictly convex surfaces $M_t$ 
in $\R^3$, flowing according to $\dt X=-H^4\nu$, 
$$\max\limits_{M_t}\frac{\left(\lx^2+\lx\ly+\ly^2\right)
(\lx+\ly)^6(\lx-\ly)^2}
{2\lx^2\ly^2}$$
is non-increasing in time.
\end{theorem}
\begin{proof}
We proceed as above.
\begin{align*}
F=&H^4,\displaybreak[1]\\
F^{ij}=&4H^3g^{ij},\displaybreak[1]\\
F^{ij,\,kl}=&12H^2g^{ij}g^{kl},\displaybreak[1]\\
\tilde w=&\log\left(\frac{-H^{10}+H^8\A2+2H^6\Akl22}
{H^4-2H^2\A2+\Akl22}\right)\displaybreak[1]\\
=&\log\left(\frac{\left(\lx^2+\lx\ly+\ly^2\right)(\lx+\ly)^6(\lx-\ly)^2}
{2\lx^2\ly^2}\right).
\end{align*}
In a critical point of $\tilde w$, we get
\begin{align*}
h_{22;\,1}=&a_1h_{11;\,1}\\
=&\frac{\ly}{\lx}\frac{8\lx^4+3\lx^3\ly+2\lx^2\ly^2-3\lx\ly^3+2\ly^4}
{8\ly^4+3\ly^2\lx+2\ly^2\lx^2-3\ly\lx^3+2\lx^4}h_{11;\,1},\displaybreak[1]\\
\oprokl H=&H^4\A2+12H^2(1+a_1)^2\cdot h_{11;\,1}^2
+(\ldots)\cdot h_{22;\,2}^2,
\displaybreak[1]\\
\opr{\A2}=&8H^3\Akl22-6H^4\A3\\
&-8H^3\left(1+3a_1^2\right)\cdot h_{11;\,1}^2
+24H^2\lx(1+a_1)^2\cdot h_{11;\,1}^2\\
&+(\ldots)\cdot h_{22;\,2}^2.
\end{align*}
\begin{align*}
-F^{ij}H_{;\,i}H_{;\,j}=
&-4H^3(1+a_1)^2\cdot h_{11;\,1}^2
+(\ldots)\cdot h_{22;\,2}^2,\displaybreak[1]\\
-F^{ij}\left(\left(\A2\right)_{;\,i}H_{;\,j}
+\left(\A2\right)_{;\,j}H_{;\,i}\right)=
&-16H^3(\lx+\ly a_1)(1+a_1)\cdot h_{11;\,1}^2\\
&+(\ldots)\cdot h_{22;\,2}^2,\displaybreak[1]\\
-F^{ij}\left(\A2\right)_{;\,i}\left(\A2\right)_{;\,j}=
&-16H^3(\lx+\ly a_1)^2\cdot h_{11;\,1}^2
+(\ldots)\cdot h_{22;\,2}^2,\\
\end{align*}

\begin{align*}
\oprokl{\tilde w}=&\left(\frac{-10H^9+8H^7\A2+12H^5\Akl22}
{-H^{10}+H^8\A2+2H^6\Akl22}\right.\\
&\qquad\left.-\frac{4H^3-4H\A2}{H^4-2H^2\A2+\Akl22}\right)\cdot\\
&\quad\cdot\left(\oprokl H\right)\displaybreak[1]\\
&+\left(\frac{H^8+4H^6\A2}{-H^{10}+H^8\A2+2H^6\Akl22}
-\frac{-2H^2+2\A2}{H^4-2H^2\A2+\Akl22}\right)\cdot\\
&\quad\cdot\left(\opr{\A2}\right)\displaybreak[1]\\
&+\left(\frac{-90H^8+56H^6\A2+60H^4\Akl22}
{-H^{10}+H^8\A2+2H^6\Akl22}\right.\\
&\qquad\left.-\frac{12H^2-4\A2}{H^4-2H^2\A2+\Akl22}\right)\cdot\\
&\quad\cdot\left(-F^{ij}H_{;\,i}H_{;\,j}\right)\displaybreak[1]\\
&+\left(\frac{8H^7+24H^5\A2}{-H^{10}+H^8\A2+2H^6\Akl22}
-\frac{-4H}{H^4-2H^2\A2+\Akl22}\right)\cdot\\
&\quad\cdot\left(-F^{ij}\left(H_{;\,i}\A2_{;\,j}+H_{;\,j}\A2_{;\,i}
\right)\right)\displaybreak[1]\\
&+\left(\frac{4H^6}{-H^{10}+H^8\A2+2H^6\Akl22}
-\frac{2}{H^4-2H^2\A2+\Akl22}\right)\cdot\\
&\quad\cdot\left(-F^{ij}\A2_{;\,i}\A2_{;\,j}\right)\displaybreak[1]\\
=&\frac{-3(\lx-\ly)^2(\lx+\ly)^3\lx\ly}{\lx^2+\lx\ly+\ly^2}\umbruch\\
&-\frac{4\ly(\lx+\ly)^5}{\lx^3(\lx^2+\lx\ly+\ly^2)(\lx-\ly)^2}\cdot
\displaybreak[1]\\
&\quad\cdot\frac1{(2\lx^4-3\lx^3\ly+2\lx^2\ly^2+3\lx\ly^2+8\ly^4)^2}
\cdot\displaybreak[1]\\
&\quad\cdot\left(4\lx^{10}+202\lx^9\ly-447\lx^8\ly^2+809\lx^7\ly^3
-16\lx^6\ly^4+696\lx^5\ly^5\right.\\
&\qquad\left.-511\lx^4\ly^6+161\lx^3\ly^7-78\lx^2\ly^8+4\lx\ly^9+40\ly^{10}
\right)\cdot h_{11;\,1}^2\displaybreak[1]\\
&+(\ldots)\cdot h_{22;\,2}^2\displaybreak[1]\\
\le&0.
\end{align*}
Here we used once more a computer algebra system and Sturm's theorem
to obtain the last inequality.
\end{proof}
Theorem \ref{H4 thm} follows. 


\section{$\vert A\vert^2+\beta H^2$-Flow}

\begin{theorem}\label{A2H2 thm}
A smooth closed strictly convex surface in $\R^3$, contracting
with normal velocity $\A2+\beta H^2$, $0\le\beta\le5$, 
converges to a round point in finite time.
\end{theorem}

\begin{theorem}
For a family of smooth closed strictly convex surfaces $M_t$ 
in $\R^3$, flowing according to $\dot X=-\left(\A2+\beta H^2\right)\nu$, 
$0\le\beta\le5$,
$$\max\limits_{M_t}
\frac{(\lambda_1+\lambda_2)(\lambda_1-\lambda_2)^2}{2\lambda_1\lambda_2}$$
is non-increasing in time.
\end{theorem}
\begin{proof}
Similarly as above, we obtain for $F=\A2+\beta H^2$ 
\begin{align*}
\oprokl H=&-\left(\A2\right)^2+2H\A3
+\beta\A2 H^2\umbruch\\
&+2\sum\hij+2\beta\sum\hii,\displaybreak[1]\\
\opr{\A2}=&2\A2\A3\umbruch\\
&+\beta\left(4\left(\A2\right)^2H-2H^2\A3\right)\umbruch\\
&+\beta\left(-4H\sum\hij+4\sum\hii\lambda_k\right).
\end{align*}

As in the proof of Theorem \ref{A2 mon thm}, we set
$$
\tilde w=\log H+\log\left(2\A2-H^2\right)
-\log\left(H^2-\A2\right)$$
and obtain in a critical point of $\tilde w$, where
$h_{22;\,1}=a_1h_{11;\,1}$
\begin{align*}
-F^{ij}H_{;\,i}H_{;\,j}=&-2(\lx+\beta H)(1+a_1)^2\cdot h_{11;\,1}^2
+(\ldots)\cdot h_{22;\,2}^2,\umbruch\\
-F^{ij}\left(\left(\A2\right)_{;\,i}H_{;\,j}
+\left(\A2\right)_{;\,j}H_{;\,i}\right)=&
-8(\lx+\beta H)(\lx+\ly a_1)(1+a_1)\cdot h_{11;\,1}^2\\
&+(\ldots)\cdot h_{22;\,2}^2,\umbruch\\
h_{22;\,1}=&\frac\ly\lx\frac{2\lx^2+\lx\ly+\ly^2}
{2\ly^2+\ly\lx+\lx^2}h_{11;\,1}.
\end{align*}

In a critical point of $\tilde w$, we obtain 
the evolution equation
\begin{align*}
\oprokl{\tilde w}=
&\left(\frac1H-\frac{2H}{2\A2-H^2}-\frac{2H}{H^2-\A2}\right)
\cdot\left(\oprokl H\right)\displaybreak[1]\\
&+\left(\frac{2}{2\A2-H^2}+\frac1{H^2-\A2}\right)
\cdot\left(\opr{\A2}\right)\displaybreak[1]\\
&+\left(\frac{6}{2\A2-H^2}+\frac{2}{H^2-\A2}\right)F^{ij}H_iH_j
\displaybreak[1]\\
&-\frac2{H\cdot\left(2\A2-H^2\right)}
F^{ij}\left(\left(\A2\right)_iH_j+\left(\A2\right)_jH_i\right)
\displaybreak[1]\\
=&-4\frac{K^2}H-2\beta HK\displaybreak[1]\\
&-2\frac{\left(5\lambda_1^8-4\lambda_1^7\lambda_2
+46\lambda_1^6\lambda_2^2+48\lambda_1^5\lambda_2^3
+72\lambda_1^4\lambda_2^4\right)\lambda_2}
{\left(\lambda_1+\lambda_2\right)\left(\lambda_1-\lambda_2\right)^2
\left(\lambda_1^2+\lambda_1\lambda_2+2\lambda_2^2\right)^2\lambda_1^3}
h_{11;\,1}^2\displaybreak[1]\\
&-2\frac{\left(44\lambda_1^3\lambda_2^5+34\lambda_1^2\lambda_2^6
+8\lambda_1\lambda_2^7+3\lambda_2^8\right)\lambda_2}
{\left(\lambda_1+\lambda_2\right)\left(\lambda_1-\lambda_2\right)^2
\left(\lambda_1^2+\lambda_1\lambda_2+2\lambda_2^2\right)^2\lambda_1^3}
h_{11;\,1}^2\displaybreak[1]\\
&+2\beta\frac{\left(\lambda_1^4-16\lambda_1^3\lambda_2
-6\lambda_1^2\lambda_2^2-8\lambda_1\lambda_2^3-3\lambda_2^4\right)
\left(\lambda_1+\lambda_2\right)^3\lambda_2}
{\left(\lambda_1-\lambda_2\right)^2
\left(\lambda_1^2+\lambda_1\lambda_2+2\lambda_2^2\right)^2\lambda_1^3}
h_{11;\,1}^2\displaybreak[1]\\
&+(\ldots)\cdot h_{22;\,2}^2.\\
\end{align*}

For $\beta=5$, the factor in front of $h_{11;\,1}^2$ equals
\begin{align*}
&-\frac{4\lambda_2^2}{(\lambda_1+\lambda_2)(\lambda_1-\lambda_2)^2
\left(\lambda_1^2+\lambda_1\lambda_2+2\lambda_2^2\right)^2\lambda_1^3}\cdot\\
&\quad\cdot\left(28\lambda_1^7+183\lambda_1^6\lambda_2
+334\lambda_1^5\lambda_2^2+371\lambda_1^4\lambda_2^3
+272\lambda_1^3\lambda_2^4+157\lambda_1^2\lambda_2^5
+54\lambda_1\lambda_2^6+9\lambda_2^7\right).
\end{align*}
Our claim follows. 

In order to see that the range for $\beta$
is sharp for applying the maximum principle,
we observe that the terms without derivatives of the
second fundamental form require that $\beta\ge0$. 
For $\lambda_2=1$ and $\lambda_1\to\infty$,
the factor in front of $h_{11;\,1}^2$ behaves like
$-112\lambda_1^{-3}$ for $\beta=5$ and like
$-10\lambda_1^{-2}$ for $\beta=0$, so we need the
upper bound $\beta\le 5$.
\end{proof}

\begin{lemma}
For a family of smooth closed strictly convex surfaces 
$M_t\subset\R^3$, $0\le t<T$, flowing according to
$\dt X=-F\nu$ with $F=\A2+\beta H^2$, $\beta\ge0$, a positive
lower bound on the principal curvatures, $\lx,\,\ly\ge\epsilon
>0$, is preserved during the evolution.
\end{lemma}
\begin{proof}
We proceed similarly as in Lemma \ref{H lower bound}.
Dropping the term involving second derivatives of $F$ yields
\begin{align*}
\left(\dt M_{ij}-F^{kl}M_{ij;\,kl}\right)\xi^i\xi^j\ge&
2\left(\A3+\beta H\A2\right)\epsilon
-\left(\A2+\beta H^2\right)\epsilon^2\displaybreak[1]\\
=&\epsilon^4+\epsilon\lambda^3+\epsilon\lambda^2(\lambda-\epsilon)
+\beta\left(\epsilon^4+\epsilon^2\lambda^2+2\epsilon\lambda^3\right)>0.
\end{align*}
\end{proof}

Similar calculations as before give an upper bound on the
velocity for $F=\A2+\beta H^2$, $\beta\ge0$. 

\begin{lemma}
The dual function to 
$F=\A2+\beta H^2$, $\beta\ge0$, is $\alpha$-concave
for $\alpha\le-2$.
\end{lemma}
\begin{proof}
We have 
\begin{align*}
\Phi=&\frac{2K-(1+\beta)H^2}{K^2},\displaybreak[1]\\
\Phi^{ij}=&\frac1{K^2}\left(\left(-2K+2(1+\beta)H^2\right)\tilde h^{ij}
-2(1+\beta)Hg^{ij}\right),\displaybreak[1]\\
\Phi^{ij,\,kl}=&\frac1{K^2}\left(\left(2K-4(1+\beta)H^2\right)
\tilde h^{ij}\tilde h^{kl}-\left(-2K+2(1+\beta)H^2\right)
\tilde h^{ik}\tilde h^{jl}\right.\\
&\left.+4(1+\beta)H\left(g^{ij}\tilde h^{kl}+g^{kl}\tilde h^{ij}\right)
-2(1+\beta)g^{ij}g^{kl}\right).
\end{align*}
We wish to show for $\alpha\le-2$ and for symmetric matrices 
$(\eta_{ij})$, that 
$$\Phi^{ij,\,kl}\eta_{ij}\eta_{kl}\le\frac{\alpha-1}{\alpha\Phi}
\Phi^{ij}\eta_{ij}\Phi^{kl}\eta_{kl}.$$
Terms involving $\eta_{12}^2$ have the right sign.

Consider $\alpha=-2$. Then it suffices to prove the inequality
\begin{align*}
&\left((1+\beta)H^2-2K\right)\begin{pmatrix}
6\frac{\ly^2}{\lx^2}+6\beta\frac{\ly^2}{\lx^2}+4\beta\frac{\ly}{\lx}&
2\beta\\2\beta&
6\frac{\lx^2}{\ly^2}+6\beta\frac{\lx^2}{\ly^2}+4\beta\frac{\lx}{\ly}
\end{pmatrix}\ge\\
\ge&\tfrac32\begin{pmatrix}
2(1+\beta)\frac{\ly^2}{\lx}+2\beta\ly\\
2(1+\beta)\frac{\lx^2}{\ly}+2\beta\lx
\end{pmatrix}\otimes
\begin{pmatrix}
2(1+\beta)\frac{\ly^2}{\lx}+2\beta\ly\\
2(1+\beta)\frac{\lx^2}{\ly}+2\beta\lx
\end{pmatrix}\\
\ge&\tfrac32\begin{pmatrix}
\scriptstyle\left(2(1+\beta)\frac{\ly^2}{\lx}+2\beta\ly\right)^2&
\scriptstyle\left(2(1+\beta)\frac{\lx^2}{\ly}+2\beta\lx\right)
\left(2(1+\beta)\frac{\ly^2}{\lx}+2\beta\ly\right)\\
\scriptstyle\left(2(1+\beta)\frac{\lx^2}{\ly}+2\beta\lx\right)
\left(2(1+\beta)\frac{\ly^2}{\lx}+2\beta\ly\right)&
\scriptstyle\left(2(1+\beta)\frac{\lx^2}{\ly}+2\beta\lx\right)^2
\end{pmatrix}
\end{align*}
in order to obtain $\alpha$-concavity for all
$\alpha\le-2$. This inequality is fulfilled, if
$$\left\{6\lx\ly+\beta\left(4\lx^2+12\lx\ly+4\ly^2\right)
+\beta^2\left(4\lx^2+8\lx\ly+4\ly^2\right)\right\}
\begin{pmatrix}\frac{\ly}{\lx}A&-A\\ -A&\frac{\lx}{\ly}A\end{pmatrix}$$
is positive semi-definite.
\end{proof}

We want to derive precise bounds on the principal curvatures. To 
this end, we use \eqref{diff est} and \eqref{F bound}
\begin{align*}
&(2+4\beta)\left((2+4\beta)3(T-t)\right)^{-2/3}\cdot
\left(1-c\cdot(T-t)^{1/12}\right)\umbruch\\
\le&F=\A2+\beta H^2\umbruch\\
\le&\lx^2+(\lx+\abs{\lx-\ly})^2+\beta(\lx+(\lx+\abs{\lx-\ly}))^2
\umbruch\\
\le&(2+4\beta)\lx^2+c\cdot F^{1/2}\cdot F^{1/4}+c\cdot F^{1/2}
\umbruch\\
\le&(2+4\beta)\lx^2+c\cdot(T-t)^{-2/3}\cdot
\left((T-t)^{1/6}+(T-t)^{1/3}\right).
\end{align*}
We get
$$\lx\ge((2+4\beta)3(T-t))^{-1/3}\cdot
\left(1+c\cdot(T-t)^{1/12}\right)$$
and a similar upper bound follows analogously. 

We obtain Theorem \ref{A2H2 thm}.


\section{$\tr A^3$-Flow}

\begin{theorem}\label{A3 thm}
A smooth closed strictly convex surface in $\R^3$, contracting
with normal velocity $\A3$, converges to a round point in finite time.
\end{theorem}

\begin{theorem}
For a family of smooth closed strictly convex surfaces $M_t$ 
in $\R^3$, flowing according to $\dot X=-\A3\nu$, 
$$\max\limits_{M_t}\frac{\left(3\lx^2+2\lx\ly+3\ly^2\right)(\lx-\ly)^2}
{\lx\ly}$$
is non-increasing in time.
\end{theorem}
\begin{proof}
Calculations as above yield
\def\adreizaehler{-H^4+4\Akl22}
\def\adreinenner{H^2-\A2}
\begin{align*}
\tilde w=&\log\left(\frac{\left(3\lx^2+2\lx\ly+3\ly^2\right)(\lx-\ly)^2}
{2\lx\ly}\right)\displaybreak[1]\\
=&\log\left(\frac{-H^4+4\Akl22}{H^2-\A2}\right),\displaybreak[1]\\
\oprokl H=&3\A4H-2\A3\A2+3\sum(\li+\lj)\hij\displaybreak[1]\\
=&3\A4H-2\A3\A2+3\left(2\lx+2\lx a_1^2+4\ly a_1^2\right)\cdot 
h_{11;\,1}^2\displaybreak[1]\\
&+(\ldots)\cdot h_{22;\,2}^2,
\displaybreak[1]\\
\opr{\A2}=&6\A4\A2-4\Akl32\umbruch\\
&-6\sum\lk^2\hij+6\sum\lk(\li+\lj)\hij
\displaybreak[1]\\
=&6\A4\A2-4\Akl32
-6\left(\lx^2+\lx^2a_1^2+2\ly^2a_1^2\right)\cdot h_{11;\,1}^2
\displaybreak[1]\\
&+6\left(2\lx^2+4\lx\ly a_1^2+2\ly^2a_1^2\right)\cdot h_{11;\,1}^2
+(\ldots)\cdot h_{22;\,2}^2,
\displaybreak[1]\\
-F^{ij}H_{;\,i}H_{;\,j}=&-3\lx^2(1+a_1)^2\cdot h_{11;\,1}^2
+(\ldots)\cdot h_{22;\,2}^2,\displaybreak[1]\\
-F^{ij}\left(\A2\right)_{;i}\left(\A2\right)_{;j}=&
-12\lx^2(\lx+\ly a_1)^2\cdot h_{11;\,1}^2+(\ldots)\cdot h_{22;\,2}^2,
\displaybreak[1]\\
h_{22;\,1}=&\frac\ly\lx\frac{9\lx^3+\lx^2\ly+3\lx\ly^2+3\ly^3}
{9\ly^3+\ly^2\lx+3\ly\lx^2+3\lx^3}\cdot h_{11;\,1},\displaybreak[1]\\
\oprokl{\tilde w}=&\left(\frac{-4H^3}\adreizaehler
-\frac{2H}\adreinenner\right)\cdot\left(\oprokl H\right)\umbruch\\
&+\left(\frac{8\A2}\adreizaehler
+\frac1\adreinenner\right)\cdot\\
&\quad\cdot\left(\opr{\A2}\right)\umbruch\\
&+\left(\frac{-12H^2}\adreizaehler-\frac2\adreinenner\right)\cdot
\left(-F^{ij}H_{;\,i}H_{;\,j}\right)\umbruch\\
&+\frac8\adreizaehler\left(-F^{ij}\left(\A2\right)_{;\,i}
\left(\A2\right)_{;\,j}\right)\umbruch\\
=&-\frac{2\left(\lx^4-2\lx^3\ly+18\lx^2\ly^2-2\lx\ly^3+\ly^4\right)\lx\ly}
{3\lx^2+2\lx\ly+3\ly^2}
\umbruch\\
&-\frac{6\ly}{\left(3\lx^2+2\lx\ly+3\ly^2\right)
\left(\lx-\ly\right)^2\lx^3}\cdot\\
&\quad\cdot\frac1{\left(3\lx^3+3\lx^2\ly+\lx\ly^2+9\ly^3\right)^2}\cdot\\
&\quad\cdot\left(63\lx^{12}+381\lx^{11}\ly-1389\lx^{10}\ly^2
+2883\lx^9\ly^3+36\lx^8\ly^4\right.\\
&\qquad\left.+1218\lx^7\ly^5+2294\lx^6\ly^6+582\lx^5\ly^7+855\lx^4\ly^8
\right.\\
&\qquad\left.+945\lx^3\ly^9+135\lx^2\ly^{10}+135\lx\ly^{11}
+54\ly^{12}\right)\cdot h_{11;\,1}^2\umbruch\\
&+(\ldots)\cdot h_{22;\,2}^2\umbruch\\
\le&0.
\end{align*}
\end{proof}

\begin{lemma}
For a family of smooth closed strictly convex surfaces 
$M_t\subset\R^3$, $0\le t<T$, flowing according to
$\dt X=-F\nu$ with $F=\B\alpha$, $\alpha\ge2$, a positive
lower bound on the principal curvatures, $\lx,\,\ly\ge\epsilon
>0$, is preserved during the evolution.
\end{lemma}
\begin{proof}
We proceed similarly as in Lemma \ref{H lower bound}.
Once again, the term involving second derivatives of $F$ 
is non-negative. As before, we obtain
\begin{align*}
\left(\dt M_{ij}-F^{kl}M_{ij;\,kl}\right)\xi^i\xi^j\ge&
\alpha\B{\alpha+1}\epsilon-\alpha\B\alpha\epsilon^2
-\B\alpha\epsilon^2+2\epsilon\B\alpha\epsilon\displaybreak[1]\\
=&\epsilon\left(\alpha\left(\lambda^{\alpha+1}
-\epsilon\lambda^{\alpha}\right)
+\epsilon^{\alpha+1}+\epsilon\lambda^\alpha\right)>0.
\end{align*}
\end{proof}

Similar calculations as in Lemma \ref{H upper bound}, using
$F+F^{ij}h_{ij}-\alpha F^{ij}h^k_ih_{kj}\ge0$ for 
$F=\B{\beta}$, $\beta\ge2$, yield for some $c_\beta>0$
$$\frac1{c_\beta}\left(\B\beta\right)^{1/\beta}\le
\frac1{c_\beta}H\le\frac{\B{\beta+1}}{\B\beta}\le
\frac{\beta+1}{\beta}\frac1\alpha$$
and an estimate as in Lemma \ref{H upper bound} follows.

\begin{lemma}
For $\alpha>0$, the dual function to $F=\B{\alpha}=
\lx^\alpha+\ly^\alpha$ is $-\alpha$-concave.
\end{lemma}
\begin{proof}
We set $\Phi=-\B{\alpha}\cdot K^{-\alpha}$ and have to show
that
$$\Phi^{ij,\,kl}\eta_{ij}\eta_{kl}\le
\frac{-\alpha-1}{-\alpha\Phi}\Phi^{ij}\eta_{ij}
\Phi^{kl}\eta_{kl}$$
for symmetric matrices $(\eta_{ij})$. Direct computations yield that
this inequality is equivalent to
\begin{align*}
&-\alpha(\alpha-1)K^{-\alpha}\etamatrix{\lx^{\alpha-2}}0{\ly^{\alpha-2}}
-2\alpha K^{-\alpha}
\sum\limits_{r=0}^{\alpha-2}\lx^r\ly^{\alpha-2-r}\eta_{12}^2
\displaybreak[1]\\
&+\alpha^2K^{-\alpha}\etamatrix{2\lx^{\alpha-2}}{\frac1{\lx\ly}\B{\alpha}}
{2\ly^{\alpha-2}}\displaybreak[1]\\
&-\alpha^2\B{\alpha}K^{-\alpha}\etamatrix{\frac1{\lx^2}}{\frac1{\lx\ly}}
{\frac1{\ly^2}}\displaybreak[1]\\
&-\alpha\B{\alpha}K^{-\alpha}\etamatrix{\frac1{\lx^2}}0{\frac1{\ly^2}}
-2\alpha\B{\alpha}K^{-\alpha}\frac1{\lx\ly}\eta_{12}^2\displaybreak[1]\\
\le&-\frac{\alpha(\alpha+1)K^{-\alpha}}{\B{\alpha}}
\etamatrix{\frac{\ly^{2\alpha}}{\lx^2}}{\lx^{\alpha-1}\ly^{\alpha-1}}
{\frac{\lx^{2\alpha}}{\ly^2}}.
\end{align*}
Further computations show that this is fulfilled, if
$$0\le\begin{pmatrix}\lx^{\alpha-2}\ly^\alpha
&-\lx^{\alpha-1}\ly^{\alpha-1}\\
-\lx^{\alpha-1}\ly^{\alpha-1}&
\lx^\alpha\ly^{\alpha-2}\end{pmatrix}.$$
\end{proof}

Then we proceed as before.
Similar calculations as for $F=\A2+\beta H^2$ 
give for $\alpha=3$ 
$$\lx,\,\ly\le(2(1+\alpha)\cdot(T-t))^{-\fracm1{1+\alpha}}
\cdot\left(1+c\cdot(T-t)^{\fracm14\fracm{\alpha-1}{\alpha+1}}\right)$$
and a corresponding lower bound holds.
This estimate holds also for $F=\B{\alpha}$, 
$\alpha=2,\,3,\,4,\,5,\,6$.

This finishes the proof of Theorem \ref{A3 thm}.


\section{$\tr A^{\alpha}$-Flow}

\begin{theorem}\label{Aalpha thm}
A smooth closed strictly convex surface in $\R^3$, contracting
with normal velocity $\A4$, $\A5$, or $\A6$, 
converges to a round point in finite time.
\end{theorem}

\begin{theorem}
For a family of smooth closed strictly convex surfaces $M_t$ 
in $\R^3$, flowing according to $\dot X=-\B{\alpha+2}\nu$, 
$\alpha=2,\,3,\,4$,
$$\max\limits_{M_t}
\frac{\left(\lx^\alpha+\ly^\alpha\right)(\lx+\ly)
(\lx-\ly)^2}{\lx\ly}$$
is non-increasing in time.
\end{theorem}
\begin{proof}
One might conjecture, that this quantity is also monotone
for other values of $\alpha$. For $\alpha=0$, corresponding to
$F=\A2$, we have already checked that in Theorem \ref{A2 mon thm}. 
Further computations for $\alpha=1,\,5,\,6,\,7$ suggest, however, 
that this quantity is not monotone for these values of $\alpha$.

We obtain
\begin{align*}
\oprokl H=&(\alpha+2)\B{\alpha+3}H-(\alpha+1)\B{\alpha+2}\A2
\displaybreak[1]\\
&+(\alpha+2)\sum\limits_{r=0}^\alpha\sum\limits_{i,\,j,\,k=1}^2
\lambda_i^r\lambda_j^{\alpha-r}
\hij,\displaybreak[1]\\
\opr{\A2}=&2(\alpha+2)\B{\alpha+3}\A2-2(\alpha+1)\B{\alpha+2}\A3
\displaybreak[1]\\
&-2(\alpha+2)\sum\limits_{i,\,j,\,k=1}^2\lk^{\alpha+1}\hij
\displaybreak[1]\\
&+2(\alpha+2)\sum\limits_{r=0}^\alpha
\sum\limits_{i,\,j,\,k=1}^2\lambda_i^r\lambda_j^{\alpha-r}
\hij\lk,\umbruch\\
\opr{\B\alpha}=&\alpha(\alpha+2)\B{\alpha+3}\B\alpha
-\alpha(\alpha+1)\B{\alpha+2}\B{\alpha+1}\umbruch\\
&-\alpha(\alpha+2)\sum\limits_{r=0}^{\alpha-2}\lk^{\alpha+1}
\li^{\alpha-2-r}\lj^r\hij\umbruch\\
&+\alpha(\alpha+2)\sum\limits_{r=0}^\alpha\li^r\lj^{\alpha-r}
\lk^{\alpha-1}\hij,
\end{align*}
\begin{align*}
\tilde w=&\log\left(\frac{\left(\lx^\alpha+\ly^\alpha\right)(\lx+\ly)
(\lx-\ly)^2}{\lx\ly}\right)\displaybreak[1]\\
\equiv&\log A+\log B+\log C-\log D,\displaybreak[1]\\
\oprokl {\tilde w}=&\frac1A\left(\oprokl A\right)
+\frac1B\left(\oprokl B\right)\displaybreak[1]\\
&+\frac1C\left(\oprokl C\right)
-\frac1D\left(\oprokl D\right)\displaybreak[1]\\
&+\frac2{B^2}F^{ij}B_{;\,i}B_{;\,j}+\frac2{C^2}
F^{ij}C_{;\,i}C_{;\,j}
+\frac1{BC}F^{ij}(B_{;\,i}C_{;\,j}+B_{;\,j}C_{;\,i})
\displaybreak[1]\\
&-\frac1{BD}F^{ij}(B_{;\,i}D_{;\,j}+B_{;\,j}D_{;\,i})
-\frac1{CD}F^{ij}(C_{;\,i}D_{;\,j}+C_{;\,j}D_{;\,i})\displaybreak[1]\\
=&\frac1{\tr A^\alpha}\cdot\left(\opr{\tr A^\alpha}\right)\\
&+\left(\frac1H-\frac{2H}{2\A2-H^2}-\frac{2H}{H^2-\A2}\right)\cdot
\left(\oprokl H\right)\displaybreak[1]\\
&+\left(\frac2{2\A2-H^2}+\frac1{H^2-\A2}\right)\cdot
\left(\oprokl{\A2}\right)\displaybreak[1]\\
&+\left(\frac2{2\A2-H^2}+\frac2{H^2-\A2}+\frac2{H^2}
+\frac{8H^2}{\left(2\A2-H^2\right)^2}\right.\\
&\quad\left.
-\frac{4H}{H\left(2\A2-H^2\right)}
-\frac{4H}{H\left(H^2-\A2\right)}\right.\\
&\quad\left.
+\frac{8H^2}{\left(2\A2-H^2\right)\left(H^2-\A2\right)}\right)
\cdot F^{ij}H_{;\,i}H_{;\,j}\displaybreak[1]\\
&+\left(\frac8{\left(2\A2-H^2\right)^2}+\frac4{\left(2\A2-H^2\right)
\left(H^2-\A2\right)}\right)\cdot\\
&\quad\cdot
F^{ij}\left(\A2\right)_{;\,i}\left(\A2\right)_{;\,j}\displaybreak[1]\\
&+\left(\frac1{H\left(H^2-\A2\right)}
-\frac{8H}{\left(2\A2-H^2\right)^2}
+\frac2{H\left(2\A2-H^2\right)}\right.\\
&\quad\left.
-\frac{6H}{\left(2\A2-H^2\right)\left(H^2-\A2\right)}\right)\cdot\\
&\quad\cdot F^{ij}\left(H_{;\,i}\left(\A2\right)_{;\,j}
+H_{;\,j}\left(\A2\right)_{;\,i}\right).
\end{align*}
Plugging this into a computer algebra program yields\\
for $\alpha=2$, corresponding to $F=\A4$,
\begin{align*}
h_{22;\,1}=&\frac{\ly}{\lx}
\frac{4\lx^4+\lx^3\ly+\lx^2\ly^2+\lx\ly^3+\ly^4}
{4\ly^4+\ly^3\lx+\ly^2\lx^2+\ly\lx^3+\lx^4}h_{11;\,1},\displaybreak[1]\\
\oprokl{\tilde w}=&\frac{-24\lx^4\ly^4}
{\left(\lx+\ly\right)\left(\lx^2+\ly^2\right)}\displaybreak[1]\\
&+\frac{-4\ly}{\left(\lx-\ly\right)^2\left(\lx+\ly\right)
\lx^3\left(\lx^2+\ly^2\right)}\cdot\\
&\quad\cdot\frac1
{\left(\lx^4+\lx^3\ly+\lx^2\ly^2+\lx\ly^3+4\ly^4\right)^2}\cdot\\
&\quad\cdot\left(11\lx^{16}+24\lx^{15}\ly+39\lx^{14}\ly^2
-328\lx^{13}\ly^3+482\lx^{12}\ly^4\right.\\
&\qquad\left.+192\lx^{11}\ly^5+215\lx^{10}\ly^6+236\lx^9\ly^7
+432\lx^8\ly^8+200\lx^7\ly^9\right.\\
&\qquad\left.+173\lx^6\ly^{10}+144\lx^5\ly^{11}+158\lx^4\ly^{12}
+32\lx^3\ly^{13}+21\lx^2\ly^{14}\right.\\
&\qquad\left.+12\lx\ly^{15}+5\ly^{16}\right)\cdot h_{11;\,1}^2
\displaybreak[1]\\
&+(\ldots)\cdot h_{22;\,2}^2,
\end{align*}
for $\alpha=3$ $\left(F=\A5\right)$
\begin{align*}
h_{22;\,1}=&\frac{\ly}{\lx}
\frac{5\lx^4-4\lx^3\ly+2\lx^2\ly^2+\ly^4}
{5\ly^4-4\ly^3\lx+2\ly^2\lx^2+\lx^4}h_{11;\,1},\displaybreak[1]\\
\oprokl{\tilde w}=&\frac{-4\left(2\lx^4-7\lx^3\ly+12\lx^2\ly^2
-7\lx\ly^3+2\ly^4\right)\lx^2\ly^2}{\left(\lx^2-\lx\ly+\ly^2\right)}
\displaybreak[1]\\
&+\frac{-10\ly}{\lx^3\left(\lx-\ly\right)^2
\left(\lx^4+2\lx^2\ly^2-4\lx\ly^3+5\ly^4\right)^2
\left(\lx^2-\lx\ly+\ly^2\right)}\cdot\\
&\quad\cdot\left(7\lx^{16}-65\lx^{15}\ly+397\lx^{14}\ly^2
-1295\lx^{13}\ly^3+2464\lx^{12}\ly^4\right.\\
&\qquad\left.-2981\lx^{11}\ly^5+2645\lx^{10}\ly^6-2007\lx^9\ly^7
+1510\lx^8\ly^8\right.\\
&\qquad\left.-1011\lx^7\ly^9+583\lx^6\ly^{10}-309\lx^5\ly^{11}
+176\lx^4\ly^{12}-71\lx^3\ly^{13}\right.\\
&\qquad\left.+23\lx^2\ly^{14}-5\lx\ly^{15}+3\ly^{16}\right)
\cdot h_{11;\,1}^2\displaybreak[1]\\
&+(\ldots)\cdot h_{22;\,2}^2,
\end{align*}
and for $\alpha=4$ $\left(F=\A6\right)$
\begin{align*}
h_{22;\,1}=&\frac{\ly}{\lx}
\frac{6\lx^6+\lx^5\ly-3\lx^4\ly^2+2\lx^2\ly^4+\lx\ly^5+\ly^6}
{6\ly^6+\ly^5\lx-3\ly^4\lx^2+2\ly^2\lx^4+\ly\lx^5+\lx^6}h_{11;\,1},
\displaybreak[1]\\
\oprokl{\tilde w}=&\frac{-10\left(\lx^8-2\lx^6\ly^2+6\lx^4\ly^4-2\lx^2\ly^6
+\ly^8\right)\lx^2\ly^2}{\left(\lx+\ly\right)\left(\lx^4+\ly^4\right)}
\displaybreak[1]\\
&+\frac{-6\ly}{\lx^3\left(\lx-\ly\right)^2\left(\lx+\ly\right)
\left(\lx^4+\ly^4\right)}\cdot\umbruch\\
&\quad\cdot\frac1{\left(\lx^6+\lx^5\ly+2\lx^4\ly^2-3\lx^2\ly^4
+\lx\ly^5+6\ly^6\right)^2}\cdot\umbruch\\
&\quad\cdot\left(17\lx^{24}-124\lx^{23}\ly+218\lx^{22}\ly^2
+646\lx^{21}\ly^3-642\lx^{20}\ly^4\right.\\
&\qquad\left.-2586\lx^{19}\ly^5+2536\lx^{18}\ly^6+3576\lx^{17}\ly^7
-2411\lx^{16}\ly^8\right.\\
&\qquad\left.-2928\lx^{15}\ly^9+1524\lx^{14}\ly^{10}
+1724\lx^{13}\ly^{11}+548\lx^{12}\ly^{12}\right.\\
&\qquad\left.-276\lx^{11}\ly^{13}-696\lx^{10}\ly^{14}
-8\lx^9\ly^{15}+499\lx^8\ly^{16}+236\lx^7\ly^{17}\right.\\
&\qquad\left.+146\lx^6\ly^{18}-66\lx^5\ly^{19}-2\lx^4\ly^{20}
+46\lx^3\ly^{21}+48\lx^2\ly^{22}\right.\\
&\qquad\left.+16\lx\ly^{23}+7\ly^{24}\right)\cdot h_{11;\,1}^2
\displaybreak[1]\\
&+(\ldots)\cdot h_{22;\,2}^2.
\end{align*}
In each case, Sturm's algorithm yields, that the 
right-hand side is non-positive.
\end{proof}
Theorem \ref{Aalpha thm} follows. 


\section{$H|A|^2$-Flow}

\begin{theorem}\label{A2H thm}
A smooth closed strictly convex surface in $\R^3$, contracting
with normal velocity $H\A2$, converges to a round point in finite time.
\end{theorem}

\begin{theorem}
For a family of smooth closed strictly convex surfaces $M_t$ 
in $\R^3$, flowing according to $\dt X=-H\A2\nu$, 
$$\max\limits_{M_t}\frac{(\lx+\ly)^2(\lx-\ly)^2}{2\lx\ly}$$
is non-increasing in time.
\end{theorem}
\begin{proof}
We proceed as above. 
\begin{align*}
w=&\log\left(\frac{(\lx+\ly)^2(\lx-\ly)^2}{2\lx\ly}\right)
\displaybreak[1]\\
=&\log\left(\frac{-H^4+2\A2H^2}{H^2-\A2}\right),\displaybreak[1]\\
h_{22;\,1}=&\frac{3\lx^2+\ly^2}{3\ly^2+\lx^2}\frac\ly\lx
h_{11;\,1},\displaybreak[1]\\
\oprokl H=&2H^2\A3-H\Akl22\\&+2\sum(\li+\lj)\hii
+2H\sum\hij\displaybreak[1]\\
=&2H^2\A3-H\Akl22\\
&+4(\lx+\ly a_1)(1+a_1)\cdot h_{11;\,1}^2
+2H(1+3a_1^2)\cdot h_{11;\,1}^2\\
&+(\ldots)\cdot h_{22;\,2}^2,\displaybreak[1]\\
\opr{\A2}=&2\Akl23-2\A2\sum\hij+4\sum(\li+\lj)\lk\hii\displaybreak[1]\\
=&2\Akl23-2\A2(1+3a_1^2)\cdot h_{11;\,1}^2\\&
+8\lx(\lx+\ly a_1)(1+a_1)\cdot h_{11;\,1}^2\\
&+(\ldots)\cdot h_{22;\,2}^2,\displaybreak[1]\\
-F^{ij}H_{;\,i}H_{;\,j}=&-\sum\left(\A2+2H\lk\right)\hii\displaybreak[1]\\
=&-\A2(1+a_1)^2\cdot h_{11;\,1}^2-2H\lx(1+a_1)^2\cdot h_{11;\,1}^2\\
&+(\ldots)\cdot h_{22;\,2}^2,
\end{align*}
\begin{align*}
-F^{ij}\left(\left(\A2\right)_{;\,i}H_{;\,j}
+\left(\A2\right)_{;\,j}H_{;\,i}\right)
=&-2\sum\left(\A2+2H\lk\right)(\li+\lj)\hii\\
=&-4\A2(\lx+\ly a_1)(1+a_1)\cdot h_{11;\,1}^2\\
&-8H(\lx+\ly a_1)(1+a_1)\lx
h_{11;\,1}^2\\
&+(\ldots)\cdot h_{22;\,2}^2,
\end{align*}
\begin{align*}
\oprokl w=&\left(\frac{-4H^3+4\A2H}{-H^4+2\A2H^2}-\frac{2H}{H^2-\A2}\right)
\left(\oprokl H\right)\\
&+\left(\frac{2H^2}{-H^4+2\A2H^2}+\frac1{H^2-\A2}\right)
\left(\opr{\A2}\right)\\
&+\left(\frac{-12H^2+4\A2}{-H^4+2\A2H^2}-\frac2{H^2-\A2}\right)
\left(-F^{ij}H_{;\,i}H_{;\,j}\right)\\
&+\frac{4H}{-H^4+2\A2H^2}
\left(-F^{ij}\left(\left(\A2\right)_{;\,i}H_{;\,j}
+\left(\A2\right)_{;\,j}H_{;\,i}\right)\right)\displaybreak[1]\\
=&-8\lx^2\ly^2\\
&-\frac{4\ly}{\left(\lx^2+3\ly^2\right)^2(\lx-\ly)^2\lx^3}\cdot\\
&\qquad\cdot\left(6\lx^8-39\lx^7\ly+91\lx^6\ly^2+\lx^5\ly^3
+91\lx^4\ly^4+3\lx^3\ly^5\right.\\
&\qquad\quad\left.+33\lx^2\ly^6+3\lx\ly^7+3\ly^8\right)
\cdot h_{11;\,1}^2\displaybreak[1]\\
&+(\ldots)\cdot h_{22;\,2}^2\umbruch\\
\le&0.
\end{align*}
\end{proof}

\begin{lemma}
For a family of smooth closed strictly convex surfaces 
$M_t\subset\R^3$, $0\le t<T$, flowing according to
$\dt X=-F\nu$ with $F=H\A2$, a positive
lower bound on the principal curvatures, $\lx,\,\ly\ge\epsilon
>0$, is preserved during the evolution.
\end{lemma}
\begin{proof}
We proceed similarly as in Lemma \ref{H lower bound}
and compute
\begin{align*}
F^{kl,\,rs}h_{kl;\,1}h_{rs;\,1}=&
(6\lx+2\ly)\cdot h_{11;\,1}^2+4(\lx+\ly)\cdot h_{11;\,1}h_{22;\,1}\\
&+4(\lx+\ly)\cdot h_{11;\,2}^2+(2\lx+6\ly)\cdot h_{22;\,1}^2\ge0,\umbruch\\
\left(\dt M_{ij}-F^{kl}M_{ij;\,kl}\right)\xi^i\xi^j\ge&
\left(\Akl22+2H\A3\right)\cdot h_{ij}\xi^i\xi^j\\
&-4H\A2h^k_ih_{kj}\xi^i\xi^j+2\epsilon H\A2h_{ij}\xi^i\xi^j
\displaybreak[1]\\
=&\epsilon^5+3\epsilon\lambda^4>0.
\end{align*}
\end{proof}

Similar calculations as in Lemma \ref{H upper bound} using
$F+F^{ij}h_{ij}-\alpha F^{ij}h^k_ih_{kj}\ge0$ for 
$F=H\A2$ give
\begin{align*}
\alpha\left(\Akl22+2H\A3\right)\le&4H\A2,\displaybreak[1]\\
\tfrac1c\left(H\A2\right)^{1/3}\le\tfrac12H
\le\frac{\A2}H\le&\frac4\alpha
\end{align*}
and an estimate as in Lemma \ref{H upper bound} follows.

\begin{lemma}
The dual function to $F=H\A2$ is $-3$-concave.
\end{lemma}
\begin{proof}
We set $\Phi=-H\A2 K^{-3}$ and want to prove that
$$\Phi^{ij,\,kl}\eta_{ij}\eta_{kl}\le\frac4{3\Phi}\Phi^{ij}\eta_{ij}
\Phi^{kl}\eta_{kl}.$$
We compute
\begin{align*}
\Phi^{ij}=&-\A2 K^{-3}g^{ij}-2H K^{-3}h^{ij}+3H\A2K^{-3}\tilde h^{ij},
\displaybreak[1]\\
\Phi^{ij,\,kl}=&-2K^{-3}\left(g^{ij}h^{kl}+h^{ij}g^{kl}\right)
+3\A2K^{-3}\left(g^{ij}\tilde h^{kl}+\tilde h^{ij}g^{kl}\right)\\
&+6HK^{-3}\left(h^{ij}\tilde h^{kl}+\tilde h^{ij}h^{kl}\right)\\
&-2HK^{-3}g^{ik}g^{jl}-9H\A2 K^{-3}\tilde h^{ij}\tilde h^{kl}
-3H\A2 K^{-3}\tilde h^{ik}\tilde h^{jl}.
\end{align*}
We have to check that
\begin{align*}
&-2H\A2\symmatrix{2\lx}{\lx+\ly}{2\ly}
+3H\Akl22\symmatrix{\frac2\lx}{\frac1\lx+\frac1\ly}{\frac2\ly}\\
&+6H^2\A2\symmatrix{2}{\frac\lx\ly+\frac\ly\lx}{2}
-2H^2\A2\symmatrix101\\
&-9H^2\Akl22\symmatrix{\frac1{\lx^2}}{\frac1{\lx\ly}}{\frac1{\ly^2}}
-3H^2\Akl22\symmatrix{\frac1{\lx^2}}0{\frac1{\ly^2}}\\
\le&-\frac43\symmatrix
{\scriptstyle\left(2\ly^2+\lx\ly+\frac{3\ly^3}\lx\right)^2}
{\scriptstyle\left(2\ly^2+\lx\ly+\frac{3\ly^3}\lx\right)\cdot
\left(2\lx^2+\lx\ly+\frac{3\lx^3}\ly\right)}
{\scriptstyle\left(2\lx^2+\lx\ly+\frac{3\lx^3}\ly\right)^2}.
\end{align*}
This is equivalent to
$$0\le\left(2\lx^4+\tfrac{20}3\lx^3\ly+\tfrac{44}3\lx^2\ly^2
+\tfrac{20}3\lx\ly^3+2\ly^4\right)
\symmatrix{\frac\ly\lx}{-1}{\frac\lx\ly}.$$
\end{proof}

Calculations as before show that
$$(16(T-t))^{-1/4}\cdot\left(1-c\cdot(T-t)^{1/8}\right)
\le\lx,\,\ly\le
(16(T-t))^{-1/4}\cdot\left(1+c\cdot(T-t)^{1/8}\right).$$
This finishes the proof of Theorem \ref{A2H thm}.


\section{$\abs A^4$-Flow}

\begin{theorem}\label{A4 thm}
A smooth closed strictly convex surface in $\R^3$, contracting
with normal velocity $\abs A^4$, 
converges to a round point in finite time.
\end{theorem}

\begin{theorem}
For a family of smooth closed strictly convex surfaces $M_t$ 
in $\R^3$, flowing according to $\dt X=-\abs A^4\nu=-\Akl22\nu$, 
$$\max\limits_{M_t}\frac{\left(\lx^4+2\lx^3\ly+4\lx^2\ly^2
+2\lx\ly^3+\ly^4\right)(\lx-\ly)^2}{(\lx+\ly)\lx\ly}$$
is non-increasing in time.
\end{theorem}
\begin{proof}
We calculate
\def\zaehlera22{{-H^6+2\A2 H^4-\Akl22 H^2+2\Akl23}}
\def\nennera22{{H^3-\A2 H}}
\begin{align*}
\tilde w=&\log\left(\frac{\left(\lx^4+2\lx^3\ly+4\lx^2\ly^2
+2\lx\ly^3+\ly^4\right)(\lx-\ly)^2}{(\lx+\ly)\lx\ly}\right)
\displaybreak[1]\\
=&\log\left(\frac{\zaehlera22}{\nennera22}\right),\displaybreak[1]\\
\oprokl H=&4H\A2\A3-3\Akl23\displaybreak[1]\\
&+8\sum\li\lj\hii+4\A2\sum\hij,\displaybreak[1]\\
\opr{\A2}=&2\Akl22\A3+16\sum\li\lj\lk\hii,
\end{align*}
\begin{align*}
\oprokl {\tilde w}=&\left(\frac{-6H^5+8\A2 H^3-2\Akl22 H}{\zaehlera22}
-\frac{3H^2-\A2}{\nennera22}\right)\cdot\\
&\quad\cdot\left(\oprokl H\right)\displaybreak[1]\\
&+\left(\frac{2H^4-2\A2 H^2+6\Akl22}{\zaehlera22}
-\frac{-H}{\nennera22}\right)\cdot\\
&\quad\cdot\left(\opr{\A2}\right)\displaybreak[1]\\
&+\left(\frac{-30H^4+24\A2 H^2-2\Akl22}{\zaehlera22}
-\frac{6H}{\nennera22}\right)\cdot\\
&\quad\left(-F^{ij}H_{;\,i}H_{;\,j}\right)\displaybreak[1]\\
&+\frac{-2H^2+12\A2}{\zaehlera22}\cdot\\
&\quad\cdot
\left(-F^{ij}\left(\A2\right)_{;\,i}\left(\A2\right)_{;\,j}\right)
\displaybreak[1]\\
&+\left(\frac{8H^3-4\A2 H}{\zaehlera22}
-\frac{-1}{\nennera22}\right)\cdot\\
&\quad\cdot\left(-F^{ij}\left(H_{;\,i}\left(\A2\right)_{;\,j}
+H_{;\,j}\left(\A2\right)_{;\,i}\right)\right).
\end{align*}
We use a computer algebra program and obtain
\begin{align*}
h_{22;\,1}=&\frac{\ly}{\lx}\frac
{4\lx^6+9\lx^5\ly+11\lx^4\ly^2+10\lx^3\ly^3+2\lx^2\ly^4+3\lx\ly^5+\ly^6}
{4\ly^6+9\lx\ly^5+11\lx^2\ly^4+10\lx^3\ly^3+2\lx^4\ly^2+3\lx^5\ly+\lx^6}
\cdot h_{11;\,1},\displaybreak[1]\\
\oprokl {\tilde w}=&\frac{-12
\left(3\lx^2+4\lx\ly+3\ly^2\right)\left(\lx^2+\ly^2\right)
\lx^3\ly^3}{\left(\lx+\ly\right)
\left(\lx^4+2\lx^3\ly+4\lx^2\ly^2+2\lx\ly^3+\ly^4\right)}
\displaybreak[1]\\
&+\frac{-4\ly}{\left(\lx^4+2\lx^3\ly+4\lx^2\ly^2+2\lx\ly^3+\ly^4\right)
\left(\lx-\ly\right)^2(\lx+\ly)\lx^3}\cdot\displaybreak[1]\\
&\quad\cdot\frac1{\left(\lx^6+3\lx^5\ly+2\lx^4\ly^2+10\lx^3\ly^3
+11\lx^2\ly^4+9\lx\ly^5+4\ly^6\right)^2}\cdot\displaybreak[1]\\
&\quad\cdot\left(11\lx^{22}+90\lx^{21}\ly-113\lx^{20}\ly^2
-840\lx^{19}\ly^3-1507\lx^{18}\ly^4\right.\\
&\qquad\left.+66\lx^{17}\ly^5+7465\lx^{16}\ly^6
+23136\lx^{15}\ly^7+45494\lx^{14}\ly^8\right.\\
&\qquad\left.+70100\lx^{13}\ly^9+84982\lx^{12}\ly^{10}
+85120\lx^{11}\ly^{11}+70882\lx^{10}\ly^{12}\right.\\
&\qquad\left.+52148\lx^9\ly^{13}+33938\lx^8\ly^{14}
+20928\lx^7\ly^{15}+11263\lx^6\ly^{16}\right.\\
&\qquad\left.+5490\lx^5\ly^{17}+2363\lx^4\ly^{18}+744\lx^3\ly^{19}
+193\lx^2\ly^{20}\right.\\
&\qquad\left.+42\lx\ly^{21}+5\ly^{22}\right)
\cdot h_{11;\,1}^2\displaybreak[1]\\
&+(\ldots)\cdot h_{22;\,2}^2.
\end{align*}
We apply Sturm's algorithm and obtain that the right-hand side
is non-positive.
\end{proof}

\begin{lemma}
For a family of smooth closed strictly convex surfaces 
$M_t\subset\R^3$, $0\le t<T$, flowing according to
$\dt X=-F\nu$ with $F=\abs A^4$, a positive
lower bound on the principal curvatures, $\lx,\,\ly\ge\epsilon
>0$, is preserved during the evolution.
\end{lemma}
\begin{proof}
The term involving second derivatives of $F$ is non-negative,
so we have
\begin{align*}
\left(\dt M_{ij}-F^{kl}M_{ij;\,kl}\right)\xi^i\xi^j\ge&
\A2\epsilon\left(4\A3-3\epsilon\A2\right)\displaybreak[1]\\
=&\A2\epsilon\left(\lambda^3+3\lambda^2(\lambda-\epsilon)
+\epsilon^3\right)>0.
\end{align*}
\end{proof}

Similar calculations as in Lemma \ref{H upper bound} using
$F+F^{ij}h_{ij}-\alpha F^{ij}h^k_ih_{kj}\ge0$ for 
$F=\abs A^4$ give
\begin{align*}
5\A2-4\alpha\A3\ge&0,\displaybreak[1]\\
\tfrac1c\left(\abs A^4\right)^{1/4}\le\frac{\A3}{\A2}\le\frac5{4\alpha}
\end{align*}
and an estimate as in Lemma \ref{H upper bound} follows.

The following lemma implies that the dual function to 
$F=\abs A^4$ is $-4$-concave.
\begin{lemma}
If the dual function to $F$ is $\alpha$-concave for
some $\alpha<0$, then the dual function to $F^\beta$
is $\alpha\cdot\beta$-concave for $\beta>0$.
\end{lemma}
\begin{proof}
We use similar notation as before. Assume that 
$$-G^{ij,\,kl}\eta_{ij}\eta_{kl}\le\frac{\alpha-1}{\alpha(-G)}
G^{ij}\eta_{ij}G^{kl}\eta_{kl}.$$
We have to show for $\Phi=-G^\beta$
$$\Phi^{ij,\,kl}\eta_{ij}\eta_{kl}\le\frac{\alpha\beta-1}
{\alpha\beta\Phi}\Phi^{ij}\eta_{ij}\Phi^{kl}\eta_{kl}.$$
Direct calculations yield
\begin{align*}
\Phi^{ij}=&-\beta G^{\beta-1}G^{ij},\displaybreak[1]\\
\Phi^{ij,\,kl}=&-\beta G^{\beta-1}G^{ij,\,kl}
-\beta(\beta-1)G^{\beta-2}G^{ij}G^{kl},\displaybreak[1]\\
\Phi^{ij,\,kl}\eta_{ij}\eta_{kl}=&
-\beta G^{\beta-1}G^{ij,\,kl}\eta_{ij}\eta_{kl}
-\beta(\beta-1)G^{\beta-2}\left(G^{ij}\eta_{ij}\right)^2
\displaybreak[1]\\
\le&-\beta\frac{\alpha-1}{\alpha}G^{\beta-2}
\left(G^{ij}\eta_{ij}\right)^2
-\beta(\beta-1)G^{\beta-2}\left(G^{ij}\eta_{ij}\right)^2
\displaybreak[1]\\
=&-\frac{\alpha\beta-1}{\alpha\beta G^\beta}\beta^2
G^{2\beta-2}\left(G^{ij}\eta_{ij}\right)^2\displaybreak[1]\\
=&\frac{\alpha\beta-1}{\alpha\beta\Phi}
\left(\Phi^{ij}\eta_{ij}\right)^2.
\end{align*}
\end{proof}

Calculations as before show that
$$\lx,\,\ly\le
(20(T-t))^{-4/5}\cdot\left(1+c\cdot(T-t)^{3/{20}}\right).$$
A corresponding lower estimate is also true.
Theorem \ref{A4 thm} follows.


\section{Convergence Rate}\label{conv rate}

In order to find out what the optimal convergence rate might
be, we consider the evolution eqution 
$$\dt X=-\A2\nu+2X,$$
that appropriately rescaled solutions of $\dt X=-\A2\nu$ fulfill.
As in \cite[Appendix]{OSKnutCrelle}, we represent the surfaces $M_t$ 
as graphs over the sphere with embeddings 
$$\S^2\ni x\mapsto x\cdot u(x),$$
where $u:\S^2\to\R_+$. Let $(\sigma_{ij})$ be the standard metric
on the sphere and $\left(\sigma^{ij}\right)$ its inverse. 
Then we get as in \cite{OSKnutCrelle}, 
using indices to denote covariant
derivatives on $\S^2$,
\begin{align*}
g_{ij}=&u^2(\sigma_{ij}+\phi_i\phi_j),\quad\text{where }\phi=\log u,
\umbruch\\
h_{ij}=&\frac1{uw}g_{ij}-\frac uw\phi_{ij},\quad\text{where }
w=\sqrt{1+\phi_i\sigma^{ij}\phi_j},\umbruch\\
\fracp ut=&-\A2 w+2u.
\end{align*}
We linearize our equation around the stationary solution 
$u=1$ and take $u=1+\epsilon v$. Then
\begin{align*}
\left.\frac d{d\epsilon} w\right|_{\epsilon=0}=&0,\umbruch\\
\left.\frac d{d\epsilon} g_{ij}\right|_{\epsilon=0}=&2v\sigma_{ij},\umbruch\\
\left.\frac d{d\epsilon} h_{ij}\right|_{\epsilon=0}=&v\sigma_{ij}-v_{ij}.
\end{align*}
So the linearized equation becomes 
$$\fracp vt=2\Delta v+6v.$$
There are spherical harmonics $u$ solving $\Delta u=-l(l+1)u$
for $l\in\N$ \cite{JonesSphericalHarmonics}. 
We do not need to consider $l=0$ as a corresponding
eigenfunction is positive (or negative) everywhere. Thus the 
resulting surface does not contract to a point for $t\uparrow T$
in the unrescaled setting. Similarly, we can exclude $l=1$ as
the respective eigenfunctions correspond to translations and
translated surfaces converge to infinity in the rescaled setting.

Using the ansatz eigenfunction 
multiplied with $e^{-\lambda t}$, we get $\lambda=2l(l+1)-6$.
This is positive for $l=2$, $\lambda=6$. So we should not
expect a convergence rate better than $||X|-1|\le c\cdot e^{-6t}$ after 
rescaling, corresponding to estimates like
$$r_+\le(6(T-t))^{1/3}\cdot(1+c\cdot(T-t)).$$
Here, we have assumed that $q(t)=0$.

We can still improve our convergence rate by considering
$$\max_{M_t}\frac{\left(\lx^2+\ly^2\right)^2(\lx-\ly)^2}{4\lx^2\ly^2}
=\max_{M_t}\frac{\left(2\A2-H^2\right)\Akl22}{\left(H^2-\A2\right)^2}
\equiv\max_{M_t}w.$$
Under the flow equation $\dt X=-\A2\nu$, we obtain, using 
calculations as above, that $w$ fulfills in a 
critical point of $w$
\begin{align*}
\oprokl w=&
-\frac{\left(\lx^2+\ly^2\right)^2}
{\left(\lx^3-\lx\ly^2+2\ly^3\right)^2\lx^5\ly}
\cdot\left(5\lx^9-3\lx^8\ly-6\lx^7\ly^2+26\lx^6\ly^3\right.\\
&\qquad\left.-20\lx^5\ly^4
+8\lx^4\ly^5+6\lx^3\ly^6-2\lx^2\ly^7-\lx\ly^8+3\ly^9\right)
\cdot h_{11;\,1}^2\\
&+(\ldots)\cdot h_{22;\,2}^2\umbruch\\
\le&0. 
\end{align*}
This implies convergence rates like
$$r_+\le(6(T-t))^{1/3}\cdot\left(1+c\cdot(T-t)^{1/3}\right)$$
similar to \cite{AndrewsStones}, where the author also has a
scaling invariant upper bound for $|\lx-\ly|$.
As there are no negative constant terms left for this choice
of $w$, it might be, that there is no monotone quantity as
studied in this paper, that allows to improve this convergence 
rate.


\bibliographystyle{amsplain}
\def\unterstrich{\underline{\rule{1ex}{0ex}}} \def\cprime{$'$} \def\cprime{$'$}
  \def\cprime{$'$}
\providecommand{\bysame}{\leavevmode\hbox to3em{\hrulefill}\thinspace}
\providecommand{\MR}{\relax\ifhmode\unskip\space\fi MR }
\providecommand{\MRhref}[2]{%
  \href{http://www.ams.org/mathscinet-getitem?mr=#1}{#2}
}
\providecommand{\href}[2]{#2}

\end{document}